\documentclass[12pt,reqno]{amsart}

\usepackage[dvips]{graphicx} 
\usepackage[arrow,matrix,curve]{xy}

\usepackage{amssymb, latexsym, amsmath, amscd, array, 
}

\newtheorem{theorem}{Theorem}[section]

\theoremstyle{definition}

\newtheorem{example}[theorem]{Example}

\newtheorem{question}[theorem]{Question}

\numberwithin{equation}{section}
\numberwithin{figure}{section} 
\numberwithin{table}{section}

\DeclareMathOperator{\PD}{{\rm PD}}

\DeclareMathOperator{\sys}{{\rm Sys}}

\DeclareMathOperator{\area}{{\rm Area}}

\newcommand \cqfd{\unskip\kern 6pt\penalty 500
\raise -2pt\hbox{\vrule\vbox to10pt{\hrule width 4pt
\vfill\hrule}\vrule}\par}                 

\def\adots{\mathinner{\mkern2mu\raise1pt\hbox{.}
\mkern3mu\raise4pt\hbox{.}\mkern1mu\raise7pt\hbox{.}}}
\def\hfl#1{\frac{\buildrel{#1}}{{\hbox to 12mm{\rightarrowfill}}}}
  
\def \\R^n \times \R^n
\rightarrow \R{\mathop{\R^n \times \R^n
\rightarrow \R}}

\newcommand\N {{\mathbb N}}

\newcommand\Q{{\mathbb {Q}}} \newcommand\R {{\mathbb R}}
\newcommand\RR {{\mathbb R}} 
\newcommand\T {{\mathbb T}} \newcommand\Z {{\mathbb Z}}  

\newcommand{\smat}[4]
{{\(\!\!\begin{array}{cc}{#1}\!&\!{#2}\\begin{equation*}-0.1cm]{#3}\!&\!{#4}\end{array}\!\!\)}}

\long\def\forget#1\forgotten{} %

\long\def\forgett#1\forgottent{} %

\def\circ{\mathchoice%
 {\mathrel{\raise 1pt\hbox{$\scriptstyle\mathchar"020E$}}}
 {\mathrel{\raise 1pt\hbox{$\scriptstyle\mathchar"020E$}}}
 {\mathrel{\raise 1pt\hbox{$\scriptscriptstyle\mathchar"020E$}}}
 {}
}

\newcommand{\nc}{\newcommand} \nc{\on}{\operatorname}

\nc{\df}{\on{\it df}}

\nc{\conf}{\on{conf}}

\nc{\spt}{\on{spt}}
\nc{\norm}[1]{\| #1 \|}
\nc{\parallelleer}{\norm{\ }} 
\nc{\parallelh}{\norm h} 
\nc{\parallelk}{\norm k} 
\nc{\parallelx}{\norm x} 
\nc{\parallelhrr}{\norm {h_\RR}} 
\nc{\parallelom}{\norm \omega} 
\nc{\parallelomij}{\norm {\omega_{i_j}}} 
\nc{\parallelomx}{\norm {\omega_{x}}} 
\nc{\parallelpi}{\norm \pi} 
\nc{\parallelalf}{\norm \alpha} 
\nc{\parallelalfs}{\norm {\alpha_s}} 
\nc{\parallelalfi}{\norm {\alpha_i}} 
\nc{\parallelalfij}{\norm {\alpha_{i_j}}} 
\nc{\parallelbeta}{\norm \beta} 
\nc{\parallelbetat}{\norm {\beta_t}} 
\nc{\parallelhcapalf}{\norm {h \cap \alpha}} 
\nc{\parallelPDralf}{\norm {\PD_\RR(\alpha)}} 
\nc{\strichleer}{| \  |}
\nc{\NN}{\mathbb N}

\nc{\rr}{\mbox{$\scriptstyle\mathbb R$}}

\nc{\dF}{{\it dF}} 

\nc{\DF}{{\it DF}} 

\nc{\ds}{{\it ds}} 

\nc{\dvol}{{\it dvol}}

\nc{\grad}{{\rm grad}} 
\nc{\strichw}{\|\omega\|} 
\nc{\strichwx}{|\omega_x|}
\nc{\Hess}{{\rm Hess}}

\newcommand\RRR{\mbox{I\!I\!R}}

\newcommand\RRRhuge{{\Huge\mbox{I\!I\!R}}}

\newcommand\Los{{\L}o{\'s}}

\begin{document}

\title[Meaning in classical mathematics and Intuitionism] {Meaning in
classical mathematics: is it at odds with Intuitionism?}


\author{Karin Usadi Katz}

\author[M.~Katz]{Mikhail G. Katz$^{*}$}

\address{Department of Mathematics, Bar Ilan University, Ramat Gan
52900 Israel} \email{katzmik ``at'' macs.biu.ac.il}

\thanks{$^{*}$Supported by the Israel Science Foundation grant
1294/06}



\date{\today}

\begin{abstract}
We examine the classical/intuitionist divide, and how it reflects on
modern theories of infinitesimals.  When leading intuitionist Heyting
announced that ``the creation of non-standard analysis is a standard
model of important mathematical research'', he was fully aware that he
was breaking ranks with Brouwer.  Was Errett Bishop faithful to either
Kronecker or Brouwer?  Through a comparative textual analysis of three
of Bishop's texts, we analyze the ideological and/or pedagogical
nature of his objections to infinitesimals \`a la Robinson.  Bishop's
famous ``debasement'' comment at the 1974 Boston workshop, published
as part of his {\em Crisis\/} lecture, in reality was never uttered in
front of an audience.  We compare the realist and the anti-realist
intuitionist narratives, and analyze the views of Dummett, Pourciau,
Richman, Shapiro, and Tennant.  Variational principles are important
physical applications, currently lacking a constructive framework.  We
examine the case of the Hawking--Penrose singularity theorem, already
analyzed by Hellman in the context of the Quine-Putnam
indispensability thesis.

\medskip\noindent {\large R\'esum\'e:} Nous analysons le clivage
classique/intuitionniste, et la fa\c con dont il est refl\'et\'e dans
les th\'eories modernes des infiniment petits.  Losque
l'intuitionniste Heyting annon\c ca que ``la cr\'eation de l'analyse
non-standard est un mod\`ele standard de recherche math\'ematique
importante'', il fut pleinement conscient du fait qu'il rompait les
rangs avec Brouwer.  Errett Bishop fut-il fid\`ele soit \`a Kronecker
soit \`a Brouwer?  Par le biais d'une \'etude comparative de trois
textes de Bishop, nous analysons la nature id\'eologique et/ou
p\'edagogique de ses objections aux infinit\'esimaux \`a la Robinson.
La c\'el\`ebre remarque de Bishop concernant le ``debasement'', au
colloque de Boston en 1974, et publi\'ee dans le cadre de son texte
{\em La Crise etc.\/}, en r\'ealit\'e ne fut jamais \'enonc\'ee devant
un public.  Nous comparons les narratives realiste et anti-realiste
intuitionnistes, et analysons les opinions de Dummett, Pourciau,
Richman, Shapiro, et Tennant.  Les principes variationnels constituent
des applications physiques importantes, et qui manquent actuellement
de cadre constructif.  Nous examinons le cas des th\'eor\`emes sur les
singularit\'es (de Hawking et Penrose), d\'ej\`a analys\'e par Hellman
dans le contexte de la th\`ese de l'indispensabilit\'e de
Quine-Putnam.
\end{abstract}

\maketitle

\tableofcontents

\section{Introduction}

\subsection{Cognitive bias}

The study of cognitive bias has its historical roots in a
classification proposed by Francis Bacon (1561 - 1626) of what he
called the {\em idola\/} (a Latin plural) of several kinds.  He
described these as things which obstructed the path of correct
scientific reasoning.  Of particular interest to us are his {\em Idola
fori\/} (Illusions of the Marketplace): due to confusions in the use
of language and taking some words in science to have a different
meaning than their common usage; and {\em Idola theatri\/} (Illusions
of the Theater): the following of academic dogma and not asking
questions about the world (see Bacon~\cite{Bac}).

{\em Completeness, continuity, continuum, Dedekind ``gaps''\/}: these
are terms whose common meaning is frequently conflated with their
technical meaning.  In our experience, explaining
infinitesimal-enriched extensions of the reals to a professional
mathematician typically elicits a puzzled reaction of the following
sort: ``But how can you {\em extend\/} the real numbers?  Aren't they
already {\em complete\/} by virtue of having filled in all the {\em
gaps\/} already?''  The academic dogma concerning the uniqueness of
the {\em continuum\/} defined necessarily as the complete Archimedean
ordered field in the context of ZFC has recently been challenged in
the pages of {\em Intellectica\/}, where numerous conceptions of the
continuum were presented, ranging from S.~Feferman's predicative
conception of the continuum \cite{Fe09} to J.~Bell's conception in
terms of an intuitionistic topos \cite{Be09b} following Lawvere
\cite{Law}.

To illustrate the variety of possible conceptions of the continuum,
note that traditionally, mathematicians have considered at least two
different types of continua.  These are Archimedean continua, or
A-continua for short, and infinitesimal-enriched (Bernoulli) continua,
or B-continua for short.  Neither an A-continuum nor a B-continuum
corresponds to a unique mathematical structure (see
Table~\ref{continuity}).  Thus, we have at least two distinct
implementations of an A-continuum:
\begin{itemize}
\item
the real numbers in the context of classical logic (incorporating the
law of excluded middle);
\item
Brouwer's real continuum built from ``free-choice sequences'', in the
context of intuitionistic logic.
\end{itemize}

\renewcommand{\arraystretch}{1.3}
\begin{table}
\[
\begin{tabular}[t]
{ | p{.9in} || p{1in} | p{.9in} | p{.5in} | p{.75in} | p{.5in} |}
\hline & Archimedean & Bernoullian \\ \hline\hline classical &
Stevin's \mbox{continuum}$^{\ref{simon}}$ & Robinson's
continuum$^{\ref{f4bis}}$ \\ \hline intuitionistic & Brouwer's
continuum & Lawvere's continuum$^{\ref{f5}}$ \\ \hline
\end{tabular}
\]
\caption{\textsf{Varieties of continua, mapped out according to a pair
of binary parameters: classical/intuitionistic and
Archimedean/Bernoullian}.}
\label{continuity}
\end{table}
\renewcommand{\arraystretch}{1}

John L. Bell describes a distinction within the class of
infinitesimal-enriched B-continua, in the following terms.
Historically, there were two main approaches to such an enriched
continuum, one by Leibniz, and one by B.~Nieuwentijdt.  The latter
favored nilpotent (nilsquare) infinitesimals.  J.~Bell notes:
\begin{quote}
Leibnizian infinitesimals (differentials) are realized in
[A.~Robinson's] nonstandard analysis,%
\footnote{\label{f4bis}More precisely, the Hewitt-{\L}o{\'s}-Robinson
continuum; see footnote~\ref{f7}.}
and nilsquare infinitesimals in [Lawvere's] smooth infinitesimal
analysis~\cite{Bel08, Bel09}.
\end{quote}
The latter theory relies on intuitionistic logic.%
\footnote{\label{f5}Lawvere's infinitesimals rely on intuitionistic
logic.  Thus, Lawvere's system is incompatible with classical logic,
see J.~Bell \cite{Bel09}.  The situation with Bishop's constructive
mathematics is more complex; see Subsection~\ref{evt}.}
A more recent implementation of an infinitesimal-enriched continuum
has been developed by P.~Giordano, see \cite{Gio09, Gio10b, Gio10a},
combining elements of both a classsical and an intuitionistic
continuum.  The familiar real continuum is rooted in the 16th century
work of Simon Stevin.%
\footnote{\label{simon}See A.~Malet~\cite{Mal06}, J.~Naets~\cite{Nae},
Katz \& Katz \cite{KK11c}.}

\subsection{Birkhoff's {\em opposing schools}}

J. Dauben reports that Garrett Birkhoff made the following statement
at a {\em Workshop on the evolution of modern mathematics\/} in~1974:
\begin{quote}
During the past twenty years, significant contributions to the
foundations of mathematics have been made by two opposing schools
\cite[p.~132]{Da96}.
\end{quote}
Birkhoff proceeded to identify the {\em opposing schools\/} as Abraham
Robinson's non-standard approach, and Errett Bishop's constructive
approach, which, according to Birkhoff,
\begin{quote}
attempts to reinterpret Brouwer's ``intuitionism'' in terms of
concepts of `constructive analysis'.
\end{quote}
S.~Warschawski referred to non-standard analysis as
\begin{quote}
a theory that has been described as diametrically opposed to
constructivism \cite[p.~37]{War}.
\end{quote}
In 1999, a conference \cite{Sc} sought to reunite the {\em
antipodes\/} (another term for our pair of opposing schools).  D.~Ross
\cite[p.~494]{Ross} evokes a {\em conceptual gulf\/} between the two
approaches.

Robinson's approach to infinitesimals has a coherent methodology in
the framework of the Zermelo-Fraenkel%
\footnote{\label{zf1}See footnote~\ref{zf2} on Fraenkel's comments 
on infinitesimals.}
set theory with the axiom of choice, in the context of classical logic
(incorporating the law of excluded middle).

The situation with intuitionism/constructivism is more nebulous.
Within this school, we will distinguish two significantly different
approaches.  One of the approaches can arguably be associated with
A.~Heyting's views \cite{He, He73}, as well as with Bishop's views
\cite{Bi75} as expressed in 1975.  The other approach can be
associated with Bishop's views \cite{Bi77} as expressed in 1977, as
well as with M.~Beeson \cite{Bee} and M.~Dummett~\cite{Du}.%
\footnote{Dummett's views will be examined in
Section~\ref{insurrection}.}
A related dichotomy of {\em liberal\/} versus {\em radical\/}
constructivism was analyzed by G.~Hellman~\cite{Hel93a}.  Rather than
borrowing terms from political science, we will exploit terms that are
somewhat more self-explanatory in a mathematical context, namely the
view of constructive mathematics as a {\em companion\/} versus an {\em
alternative\/} to classical mathematics.

The infinitesimal, affectionately nicknamed the {\em cholera
bacillus\/} of mathematics by G. Cantor,%
\footnote{See Meschkowski \cite[p.~505]{Me} and Dauben
\cite[p.~353]{Da95}, \cite[p.~124]{Da96}.}
has mutated into a bone of contention between our pair of antipodes, a
bone constricting the constructivist throat of our students'
aspirations, according to some; and the elegant answer to Leibniz's
dream,%
\footnote{\label{f49} Note that Schubring \cite[p.~170, 173, 187]{Sch}
attributes the first systematic use of infinitesimals as a
foundational concept, to Johann Bernoulli (rather than Leibniz).}
according to others.

The possibility of a peaceful coexistence between our pair of
antipodes is not an idle dream, but rather a matter of historical
record.  In 1954, Robinson gave a lecture at a symposium on the
``Mathematical interpretation of formal systems'' chaired by
A. Heyting, see \cite[p.~254]{Da03}.  In 1961, Robinson made public
his new idea of non-standard models for analysis, and ``communicated
this almost immediately to A[.]~Heyting'' \cite[p.~259]{Da03}.
Robinson's first paper on the subject was subsequently published in
{\em Proceedings of the Netherlands Royal Academy of
Sciences\/}~\cite{Ro61}.  An appreciation of Heyting's work appears in
Robinson's 1968 text:
\begin{quote}
Among those who have shown a consistently constructive attitude we may
mention Heyting\dots.  \cite[p.~921]{Ro68}.
\end{quote}
In 1973, Heyting praised non-standard analysis as ``a standard model
of important mathematical research'' \cite[p.~136]{He73}.%
\footnote{Heyting's appreciation of Robinson's theory is analyzed in
more detail in Subsection~\ref{APR}.}

Meanwhile, S.~Shapiro unveils an allegedly irreducible clash between
classical mathematics, on the one hand, and intuitionism as expressed
through Heyting's semantics (see \cite[p.~60]{Sh}), on the other.%
\footnote{Shapiro's views will be analyzed in Subsection~\ref{eight}.}

From the viewpoint of classical mathematics, what is a reasonable
scope of the constructivist critique of the foundations?  We will
describe the proposed dichotomy within constructivism, and argue that
only one of the approaches would be potentially amenable to a possible
reunification (or at least a state of non-belligerence) alluded to in
the title of the conference~\cite{Sc}.

While commenting on the ``role of non-standard analysis in
mathematics'', P.~Halmos writes:
\begin{quote}
For some other[\dots mathematicians], who are against it (for instance
Errett Bishop), it's an equally emotional issue \cite[p.~204]{Ha}.
\end{quote}
What was it specifically about non-standard analysis, among other
fields of classical mathematics, that prompted Bishop's opposition?
We will analyze this issue in Subsections~\ref{eb} and \ref{whyfive}.

\renewcommand{\arraystretch}{1.3} 
\begin{table}
\[
\begin{tabular}[t]
{ | p{.9in} || p{2.7in} | p{.8in} | p{.5in} | p{.75in} | p{.5in} |}
\hline field & issues & discussed in: \\ \hline\hline Mathematics &
Arguments in favor of coherence of post-LEM numerical meaning &
Section~\ref{postlem} \\ \hline Physics & Variational principles,
Hawking--Penrose theorem, Hilbert's Lagrangian & Section~\ref{HPT} \\
\hline Philosophy & Diverging insurrectional narratives: realist and
anti-realist & Section~\ref{insurrection} \\ \hline History &
Brouwer's intellectual debt to Frege's revolution and Cantor's
revolution & Section~\ref{insurrection} \\ \hline Defection & Weyl,
Heyting, Bridges arguably side with a {\em companion\/} view &
Sections~\ref{eight}, \ref{APR}, \ref{fourteen}, \ref{nine} \\ \hline
\end{tabular}
\]
\caption{\textsf{Radical Intuitionism: key issues}}
\label{summary}
\end{table}
\renewcommand{\arraystretch}{1}

Does intuitionism have an Achilles' heel?  Without attempting to
provide a definitive answer, we list potential vulnerabilities in
Table~\ref{summary}, arranged by field.

A recent analysis of the Classical {\em vs\/} Intuitionistic divide by
D.~Westerst\aa hl posits that
\begin{quote}
The typical intuitionist takes truth to be what philosophers call an
epistemic notion: roughly, something is true if it can be proved. 
\cite[p.~177]{Wes}.
\end{quote}
Such a summary would be agreed upon by many intuitionists.  But then
Westerst\aa hl concludes:
\begin{quote}
This puts computation at center stage: (intuitionistic) proofs are
computations, or directions for finding computations
\cite[p.~177]{Wes}.
\end{quote}
Westerst\aa hl's {\em reformulation\/} of constructively meaningful
mathematics in terms of ``computations'' betrays an unmistakable
influence of Bishop's perspective, which tends to conflate meaning
with computational meaning, as we will analyze below.  Similarly, on
the same page it emerges that, to Westerst\aa hl (following Dummett),
{\em classical mathematics\/} means {\em Platonism\/}.  No wonder that
his evaluation of the possibility of reaching a mutual understanding
\begin{quote}
is mostly negative: a basic {\em asymmetry\/} as to one side's ability
to achieve understanding of what the other is up to will remain
\cite[p.~177]{Wes}.
\end{quote}
Meanwhile, our own assessment is mostly {\em positive\/}, relying on a
more nuanced presentation of both sides of what we hold to be a
bridgeable divide.

\subsection{Where Weyl, Crowe, and Heyting agree}

We conclude this introduction by quoting Hermann Weyl:

\begin{quote}
It is pretty clear that our theory of the physical world is not a
description of the phenomena as we perceive them, but [rather] is a
bold symbolic construction.  However, one may be surprised to learn
that even mathematics shares this character \cite[p.~553]{Wey51}.
\end{quote}

No single formalism, whether classical or intuitionistic, represents
``a description of the [mathematical] phenomena'', to borrow Weyl's
phrase.  This sentiment is echoed by M.~Crowe:%
\footnote{\label{crowe}See footnote~\ref{curious}.}
\begin{quote}
creative mathematicians have repeatedly encountered situations that
are not resolvable by logic, that are no less dependent on brilliance
of mind and an ability to {\bf see beyond logic}, than those that have
been faced by physicists.  One example would be the situation faced by
the pioneers of the calculus, who pushed forward despite being
surrounded by inconsistencies and counter-intuitive deductions
\cite{Cr}, \cite[p.~316]{Cr92} [emphasis added--authors].
\end{quote}

In reference to classical mathematics, a similar thought is expressed
by the intuitionist Heyting, who described it as a ``curious mixture
of formal reasoning and more or less vague intuitions''.

We analyze the sources of constructive aspirations in Kronecker's
vision in Subsection~\ref{kron}.  In Subsection~\ref{kol}, we give a
gentle introduction to the intuitionistic challenge, following
A. N. Kolmogorov.  In Subsection~\ref{two}, we give some examples of
what it would mean to eliminate the reliance on the law of excluded
middle.  In Subsection~\ref{BB}, we compare Brouwerian and Bishopian
mathematics.

In Subsection~\ref{eb} we document some specific criticisms voiced by
Bishop.  In Subsection~\ref{hal}, we examine the dual aspects of
Halmos's relation to Robinson's theory.  In Subsection~\ref{whyfive},
we analyze the nature of the challenge to Bishopian constructivism
posed by a modern theory of infinitesimals.  The story behind Bishop's
``debasement'' comment is told in Subsection~\ref{debasement}.
Subsections~\ref{hypo} and~\ref{hypo2} analyze the foundational nature
of Bishop's criticisms.

The shifting ground of constructivist definitions is reviewed in
Section~\ref{shesh}.  The meaning of ``finite'' in the context of the
Goodman-Myhill theorem is reviewed in Section~\ref{six}.
Section~\ref{evt} examines the finessing of Brouwerian counterexamples
in Bishopian constructivism.  

A dichotomy of two distinct approaches within constructivism is
proposed in Subsection~\ref{51}.  Section~\ref{postlem} introduces a
dichotomy of pre-LEM and post-LEM numerical meaning.
Subsection~\ref{53} examines a philosophical outlook that can permit
an intuitionist to accept both types of meaning.

Section~\ref{fervor} deals with the fervor of Bishopian
constructivism.  Subsection~\ref{pour} examines the views of
B.~Pourciau.  A glossary of Bishopian constructivism appears in
Subsection~\ref{glossary}, followed by an analysis of a pair of
insurrectional narratives in Subsections~\ref{insurrection} and
\ref{insurrection2}.  Section~\ref{nine} deals with the views of
constructivist F.~Richman.

Section~\ref{partof} examines the issue of whether constructive
mathematics is part of classical mathematics.  Subsection~\ref{71}
examines the logical developments, including Frege's relational logic,
in the 19th century that created the conditions that made
intuitionistic revolt possible.  We focus on Heyting's intuitionism in
Subsection~\ref{eight}, where we also comment on the views of Shapiro
and Tennant (the views of the latter are also examined in
Subsection~\ref{HPT}).  Heyting's appreciation of non-standard
analysis is examined in Subsection~\ref{APR}.

Section~\ref{fourteen} deals with a challenge to constructivism
arising from natural science.  Subsection~\ref{81} deals with Bishop's
own view of the applicability of mathematics in the natural sciences.
Section~\ref{fifteen} deals with G.~Hellman's indispensability thesis
approach.  S.P. Novikov's perspective is discussed in
Subsection~\ref{83}.  The Hawking--Penrose theorem and its
foundational ramifications are analyzed in Subsection~\ref{HPT}.

\section{A crash course in intuitionism}

We present some preliminary remarks on the intuitionist project from
Kronecker to Brouwer and beyond.  More technical aspects of
constructivism are dealt with in Section~\ref{debasing}.

\subsection{The genesis of the constructive project in Kronecker}
\label{kron}

While constructivist aspirations can be detected already in a
predilection for the rationals (at the expense of broader number
systems) on the part of mathematicians both ancient and modern, the
roots of 20th century constructivism are surely to be sought in
L.~Kronecker.  J.~Boniface and N.~Schappacher recently edited an
important unpublished essay of Kronecker's.  Here Kronecker envisions
a 3-part subdivision of mathematics as follows:
\begin{quote}
la m\'ecanique, qui op\`ere avec la notion de temps, la g\'eom\'etrie,
qui \'etudie les relations spatiales ind\'ependantes du temps, et la
math\'ematique dite pure, dans laquelle n'interviennent ni le temps ni
l'espace, et que je veux appeler `arithm\'etique'%
\footnote{These comments appear on pp.~10--11 in Kronecker's
manuscript.}
\cite[p.~211]{BS}.
\end{quote}
Thus, mathematics is subdivided by Kronecker into mechanics (which
apparently corresponds to mathematical physics), geometry, and
arithmetic.  Kronecker proceeds to state the scope of his
constructivist project of the reduction of mathematics to the
integers, in the following terms:
\begin{quote}
Il souligne (p. 15) que ``tout ce qui n'appartient pas \`a la
m\'ecanique et \`a la g\'eom\'etrie, et que je veux rassembler sous
l'intitul\'e d'arithm\'etique, devrait \^etre effectivement
arithm\'etis\'e" \cite[p.~211]{BS}.
\end{quote}
Thus, the only branch of mathematics that is to be arithmetized (i.e.,
converted into a discourse about the natural numbers), according to
Kronecker, is that which falls outside the domain of both geometry and
mathematical physics.  The implication is that geometry need not, or
could not, be so arithmetized, in Kronecker's opinion.  Kronecker
provides a reasonable limitation to the scope of his constructivist
project.  Boniface and Schappacher further note:
\begin{quote}
`Arithm\'etiser' consiste, pour Kronecker, non seulement \`a r\'eduire
les objets math\'ematiques aux entiers positifs, mais aussi \`a
n'utiliser que des m\'ethodes arithm\'etiques \cite[p.~211]{BS} \dots
Les concepts de nombre n\'egatif et de nombre fractionnaire \'etant
\'evit\'es en tant que concepts fondamentaux de l'arithm\'etique,
celui de nombre irrationnel le sera a fortiori.  L'irrationalit\'e est
un concept g\'eom\'etrique et doit, selon Kronecker, rester dans le
domaine g\'eom\'etrique \cite[p.~213]{BS}.
\end{quote}
In short, geometry (including irrational numbers) falls outside the
scope of Kronecker's constructivist project.  In his moderation of
scope, Kronecker was not followed by all of his would-be disciples.

\subsection{Kolmogorov on paradoxes}
\label{kol}

Eighty years ago, A. Kolmogorov described the emergence of {\em
both\/} Hilbert's formalism and Brouwer's intuitionism, as
\begin{quote}
a reaction against the set-theoretic conception of mathematics
\cite[p.~380]{Ko}.
\end{quote}
Kolmogorov speaks of {\em great difficulties and even
contradictions\/} (ibid. p.~382) caused by {\em set-theoretic
mathematics\/}, citing Russell's paradox as an example (ibid. p.~383).
As to Zermelo's axiom of choice, he writes:
\begin{quote}
it came to a hopeless collision with the idea that mathematic[al]
existence should be based on a[n appropriate] construction [\dots]
objects whose existence is postulated by this axiom appeared to be not
only useless but sometimes destructive to the simplicity and
[rigorousness] of crucial mathematical theories.
\end{quote}
The recent republication \cite{Ko} of Kolmogorov's text was annotated
by V.~Uspenskiy.  In his note 15 on page 387, Uspenskiy illustrates
the above comment in terms of the famous {\em Banach-Tarski paradox\/}
(see \cite{Wa}).

In short, a significant mathematical minority has felt that the
paradoxes of set theory have sapped the faith in both set-theoretic
foundations and classical logic.

How is one to view the classical axiom of choice and its
counterintuitive consequences such as the Banach-Tarski paradox?
There could be at least two distinct approaches:
\begin{enumerate}
\item
one can maintain that ``infinite sets behave differently from finite
sets'' (and in particular have properties that {\em appear to be\/}
paradoxical when compared to finite sets); alternatively,
\item
some of the basic assumptions of set-theoric foundations and classical
logic need to be re-examined.
\end{enumerate}

A classically trained mathematician typically follows the first
approach.  Meanwhile, an intuitionist would follow the second
approach, and would ask the classical mathematician: ``Certainly,
infinite sets behave differently from finite sets.  But can you
exhibit a single such paradoxical infinite `set' explicitly?  The only
way of accessing such `sets' is by invoking the classical axiom of
choice and the law of excluded middle.  Since such `sets' tend to
contradict our intuition of physical space, the reasonable alternative
is to reject both the `sets' themselves, and whatever dubious axioms
and procedures are responsible for propelling them into existence in
the first place."%
\footnote{All constructivists reject the sets involved in paradoxical
decompositions \`a la Banach-Tarski.  Infinite sets (actual infinity)
are acceptable to Bishop's followers, see footnote~\ref{kronecker},
but not to Dummett's followers, see footnote~\ref{dummett}.}

The procedures in question are detailed in the next section.

\subsection{Law of excluded middle and constructive quantifiers}
\label{two}

A gentle introduction to intuitionism would necessarily seek to
clarify the law of excluded middle by providing suitable examples (see
below).  The meta-mathematical language of mainstream mathematics
relies on what is generally referred to as {\em classical logic\/},
which we will briefly denote
\[
{\rm ClLo},
\]
incorporating the Law of Excluded Middle (LEM).  The law of excluded
middle is the key ingredient in a {\em proof by contradiction\/}, as
Example~\ref{roottwo} below illustrates.  H.~Billinge writes that
E.~Bishop
\begin{quote}
seems to take it that common sense would balk at allowing the
assertion of `there exists an~$x$ such that~$P(x)$' on the flimsy
basis of having shown that a contradiction can be derived from
assuming that there is no~$x$ such that [\dots]~$P(x)$
\cite[p.~184]{Bil03}.
\end{quote}
She concludes that
\begin{quote}
the only reason that we ordinarily accept such a derivation as
evidence for such an existence claim must be that we take it that we
are dealing with a finite domain and hence that we would be able to
find the~$x$ in question by examination [i.e., ``by inspection''].%
\footnote{Of course, no computer can examine, or inspect, a collection
containing, say,~$10^{100}$ elements.  Billinge's comments should not
be interpreted as an endorsement of a strict finitism (certainly not
espoused by Bishop), but rather as a way of introducing an
intuitionistic viewpoint to a classical reader unfamiliar with this
circle of ideas.}
\end{quote}
We will illustrate the difference between the classical and the
constructive approach by means of some typical examples.

\begin{example}
\label{roottwo}
Irrationality of~$x$ is defined constructively in terms of concrete
separation between~$x$ and each rational~$\frac{m}{n}$ (separation
being expressed in terms of the denominator~$n$).  The classical proof
of the irrationality of~$\sqrt{2}$ is a proof by contradiction.
Namely, we {\em assume\/} a hypothesized
equality~$\sqrt{2}=\frac{m}{n}$, square both sides, examine the parity
of the powers of~$2$, and arrive at a contradiction.  At this stage,
irrationality is considered to have been proved, in ClLo.%
\footnote{\label{f1}The classical proof that~$\sqrt{2}$ is {\em not
rational\/} is acceptable in intuitionistic logic.  To pass from this
to the claim of its {\em irrationality\/} as defined above, requires
LEM, or more precisely the law of
trichotomy:~$(x<0)\vee(x=0)\vee(x>0)$ (see Subsection~\ref{evt} for
more details on the latter).  Our purpose here is not to present an
exhaustive treament of the intuitionistic perspective on
irrationality, but rather to give an accessible illustration of what
might be considered a numerically meaningful argument, as opposed to
an indirect proof (by contradiction).}

Alternative, {\em numerically meaningful\/} (see glossary in
Subsection~\ref{glossary}), arguments for irrationality exist.  Thus,
without exploiting the hypothetical equality~$\sqrt{2}=\frac{m}{n}$,
one can exhibit positive lower bounds for the difference
$|\sqrt{2}-\frac{m}{n}|$ in terms of the denominator~$n$, resulting in
a constructively adequate proof of irrationality.%
\footnote{For each rational~$m/n$, the integer~$2n^2$ is divisible by
an odd power of~$2$, while~$m^2$ is divisible by an even power of~$2$.
Hence~$|2n^2-m^2|\geq 1$ (here we have applied LEM to an effectively
decidable predicate over~$\Z$, or more precisely the law of
trichotomy).  Since the decimal expansion of~$\sqrt{2}$ starts with
$1.41\ldots$, we may assume~$\frac{m}{n} \leq 1.5$.  It follows that
\[
|\sqrt{2} - \tfrac{m}{n}| = \frac{|2n^2 - m^2|}{n^2 \left(
  \sqrt{2}+\tfrac{m}{n} \right)} \geq \frac{1}{n^2
  \left(\sqrt{2}+\tfrac{m}{n}\right)} \geq \frac{1}{3n^2},
\]
yielding a numerically meaningful proof of irrationality.}
A more detailed discussion of this example, by E.~Bishop, may be found
at \cite[p.~18]{Bi85}.
\end{example}

Brouwer's intuitionism, as well as most forms of constructivism, rely
on {\em intuitionistic logic}, which we will briefly denote
\[
{\rm InLo},
\]
characterized by the rejection of LEM.  The relation can be
represented graphically as follows:%
\footnote{\label{lem}The approximate equality sign~``$\approx$''
reflects the fact that there are many subsystems of classical logic
which, when LEM is adjoined to them, yield all of classical logic.
Thus, intuitionistic logic is not uniquely characterized by the
rejection of LEM.  Still, such a rejection is the main theme of
intuitionism (see Bishop's comment on the axiom of choice in
Subsection~\ref{six}).  While it may have been more correct to
write~${\rm InLo} + {\rm LEM} = {\rm ClLo}$, we chose to emphasize the
passage from classical logic to intuitionistic logic, rather than vice
versa, keeping in mind the dominant modes of logical reasoning of the
mathematical public today.}
\begin{equation}
\label{21}
{\rm ClLo} - {\rm LEM} \approx {\rm InLo}.
\end{equation}
Similarly, constructive quantifiers are different from the classical
ones, in that their domain of discourse is reduced to constructive
entities.%
\footnote{We are mostly following Bishop's views here; for a
discussion of Richman's views, see Subsections \ref{nine}
and~\ref{81}.}
Constructivism adopts a verificational interpretation of the
quantifiers.  In this context, it is relevant to mention that the
(sub)finite axiom of choice would imply LEM, see Section~\ref{six}.
To emphasize the distinction between two types of quantifiers, one
could incorporate a mention of ClLo (respectively, InLo) as part of
the quantifier notation, as the following example illustrates.

\begin{example}
Consider the following special case of the axiom of choice.  Let
$2^\R$ be the set of sets of reals, and~$2^{2^\R}$, the set of sets of
sets of reals.  Let
\[
2_{_{\coprod}}^{2^\R}
\]
be the set of sets of mutually disjoint, non-empty sets of reals.  In
this case, the classical axiom of choice asserts that
\[
\left( \forall_{\rm ClLo}^{\phantom{ClLo}}A\in 2_{_{\coprod}}^{2^\R}
\right) \left( \exists_{\rm ClLo}^{\phantom{ClLo}} S\in 2^\R \right)
\left( \forall_{\rm ClLo}^{\phantom{ClLo}} x\in A \right) \; |x\cap
S|=1.
\]
In the case of a finite set, say
\[
\{ 0, 1 \},
\]
the corresponding formula
\begin{equation}
\label{22}
\left( \forall_{\rm ClLo}^{\phantom{ClLo}}A\in
2_{_{\coprod}}^{2^{\{0,1\}}} \right) \left( \exists_{\rm
ClLo}^{\phantom{ClLo}} S\in 2^{\{0,1\}} \right) \left( \forall_{\rm
ClLo}^{\phantom{ClLo}} x\in A \right) \;\; |x\cap S|=1
\end{equation}
can be proved by induction.  Moreover, by the Goodman-Myhill theorem
(see Section~\ref{six}), formula~\eqref{22} implies LEM.%
\footnote{See footnote~\ref{bm} for a detailed discussion of this
implication.}

An intuitionistic version is obtained by replacing the subscript ClLo
by InLo (see \eqref{21}) on all quantifiers:%
\footnote{When switching to intuitionistic quantifiers, it is
customary to replace ``non-empty sets'' by ``inhabited sets'', so as
to avoid a linguistic excluded-third pitfall, see footnote~\ref{bm}.}
\begin{equation}
\label{23}
\left(
\forall_{\rm InLo}^{\phantom{InLo}}A\in
2_{_{\coprod}}^{2^{\{0,1\}}}
\right)
\left(
\exists_{\rm InLo}^{\phantom{InLo}}
S\in 2^{\{0,1\}}
\right)
\left(
\forall_{\rm InLo}^{\phantom{InLo}} x\in A 
\right)
\;\; |x\cap S|=1 
\end{equation}
Formula~\eqref{23} applies to {\em constructively\/} finite sets but
not to {\em sub\/}finite sets.%
\footnote{Subfinite sets are discussed in more detail in
Subsection~\ref{six}.}
Namely, the domain of the quantifiers is limited to constructively
defined sets only.  In particular,~\eqref{23} does not prove LEM, and
is constructively acceptable.%
\footnote{Concerning the infinite {\em constructive\/} axiom of
choice, see Bishop's remark in Subsection~\ref{six}.}
\end{example}

The intuitionistic way of thinking represents a major paradigm shift.
Many classically-trained mathematicians find it difficult to
appreciate either its utility or the well-foundedness of its trademark
distrust of proofs by contradiction.  Hopefully the examples provided
in this section provide a start for such an appreciation.

The following argument is apparently well-known in constructive
circles.%
\footnote{Communicated by Gabriel Stolzenberg.}
It concerns the claim of the consistency of the foundations of
mathematics.  The argument in favor of the claim is indirect, relying
on a proof by contradiction.  It is surprisingly short, and goes as
follows.  Consider the claim that the foundations of mathematics are
consistent.  Suppose not.  Then there is an inconsistency in the
foundations.  We get a contradiction, proving the
claim.%
\footnote{At this point the classical mathematician typically bursts
out laughing at the proof, while the constructivist may perhaps smile
at the classical mathematician.  Hopefully we have succeeded in
conveying to our gentle reader, a sense of a constructivist reception
of a typical proof by contradiction.}
\hfill\hfill~$\qed$

\subsection{Bishop and Brouwerian countexamples}
\label{BB}

E. Bishop was a leading representative of the tendency described as
intuitionist or constructivist.%
\footnote{\label{kronecker}It should be kept in mind, however, that
Bishop rejects both Kronecker's finitism (as Bishop freely exploits
the actual infinity of the natural numbers) and Brouwer's theory of
the continuum, see main text at footnote~\ref{bugaboo}.}
Among a variety of intuitionistic schools, Bishopian constructivism
has been one of the most influential, both mathematically and
philosophically.  Bishop's perspective on the Classical/Intuitionistic
divide has influenced even would-be ``outside'' observers.  Thus,
Westerst\aa hl claims that
\begin{quote} 
[Bishop's] aim was to do constructive mathematics that looked just
like ordinary mathematics, not even apparently contradicting any
classical theorems \cite[p.~185]{Wes}.
\end{quote}
However, Brouwerian counterexamples do appear in Bishop's work.  This
could not be otherwise, since a verificational interpretation of the
quantifiers necessarily results in a clash with classical mathematics.
As a matter of presentation, the conflict with classical mathematics
had been de-emphasized by Bishop.  Bishop finesses the issue of
Brouwer's theorems (e.g., that every function is continuous) by
declaring that he will only deal with uniformly continuous functions
to begin with.  In Bishopian mathematics, a circle cannot be
decomposed into a pair of antipodal sets.  A real number~$a$
satisfying~$\neg ((a\leq 0) \vee (a \geq 0))$ yields a counterexample
to the extreme value theorem.%
\footnote{A more detailed discussion may be found in
Subsection~\ref{evt}.}

\section{A kinder, gentler constructivism}
\label{kinder}

In Section~\ref{BB}, we discussed Bishop's mathematical and
philosophical influence.  Bishop put such influence to a variety of
uses.

\subsection{The {\em debasement\/} and {\em obfuscation\/} charges}
\label{eb}

Bishop's {\em Crisis\/}%
\footnote{Apparently a take-off on H.~Weyl's own {\em Crisis\/}
text~\cite{We21}, see also end of Section~\ref{fourteen}.}
essay \cite{Bi75} is cast in the form of an imaginary dialog between
Brouwer and Hilbert.  Bishop narrates a creation story of intuitionism
in the form of such an imaginary dialog, where Brouwer completely
dominates the exchange.  Indeed, Bishop's imaginary Brouwer-Hilbert
exchange is dominated by an unspoken assumption that Brouwer is the
only one who seeks ``meaning", an assumption that his illustrious
opponent is never given a chance to challenge.  

Meanwhile, Hilbert's comments in 1919 reveal clearly his attachment to
meaning which he refers to as {\em internal necessity\/}:
\begin{quote}
We are not speaking here of arbitrariness in any sense.  Mathematics
is not like a game whose tasks are determined by arbitrarily
stipulated rules.  Rather, it is a conceptual system possessing
internal necessity that can only be so and by no means otherwise
\cite[p.~14]{Hi19}.%
\footnote{\label{corry}Cited in Corry \cite{Cor}.  See also
footnote~\ref{hilbert}.}
\end{quote}

Bishop expressed his views about Robinson's infinitesimals and their
use in teaching in a brief paragraph toward the end of his text.
Following his discussion of Hilbert's formalist program, Bishop
inserted the following text:

\begin{quote}
A more recent attempt at mathematics by formal finesse is non-standard
analysis.  I gather that it has met with some degree of success,
whether at the expense of giving significantly less meaningful proofs
I do not know.  My interest in non-standard analysis is that attempts
are being made to introduce it into calculus courses.  It is difficult
to believe that {\bf debasement of meaning\/} could be carried so far%
\footnote{The story behind the insertion of this paragraph is told in
Subsection~\ref{debasement}.}
\cite[p.~513-514]{Bi75} [emphasis added--authors].
\end{quote}
Bishop's view of the introduction of Robinson-style infinitesimals in
the classroom as no less than a {\em debasement of meaning\/} (see
glossary in Subsection~\ref{glossary}), was duly noted by J.~Dauben.%
\footnote{See~\cite{Da92, Da96}.  While there was apparently a wave of
reactions to Bishop's review at the time, surprisingly little of it
has made its way into print.  Other than Dauben's essays already
mentioned, there is Keisler's own brief measured response~\cite{Ke77},
followed in the same issue of the Notices by some comments by
V.~Komkov \cite{Kom}, see footnote~\ref{komk}.  Komkov's letter was
mentioned in a MathSciNet review by P.~Smith~\cite{Sm}.  Bishop's
review was also examined briefly by M.~Davis~\cite{Dav}, see
footnote~\ref{f4}.}
Bishop's sentiments toward calculus based on Robinson-style
infinitesimals stand in sharp contrast with those of his fellow
intuitionist A.~Heyting, who felt that
\begin{quote}
[Robinson] connected [an] extremely abstract part of model theory with
a theory apparently so far apart as the elementary calculus.  In doing
so [he] threw new light on the history of the calculus by giving a
clear sense to Leibniz's%
\footnote{But see footnote~\ref{f49} for a historical clarification by
Schubring.}
notion of infinitesimals \cite[p.~136]{He73}.
\end{quote}
Heyting's views are analyzed in more detail in Subsections~\ref{eight}
and~\ref{APR}.

Bishop anchors his foundational stance in a species of mathematical
constructivism in the following terms:
\begin{quote}
To my mind, it is a major defect of our profession that we refuse to
distinguish [\dots] between integers that are computable and those that
are not [\dots] the distinction between computable and non-computable,
or constructive and non-constructive is the source of the most famous
disputes in the philosophy of mathematics
 \cite[pp.~507-508]{Bi75}.
\end{quote}
On page 511, Bishop defines a {\em limited principle of omniscience\/}
(LPO) as the supposition that

\begin{quote}
it is possible to search ``a sequence of integers to see whether they
all vanish'',
\end{quote}
and goes on to characterize the dependence on the LPO as a procedure
both Brouwer and Bishop himself reject.%
\footnote{\label{ff4}S.~Feferman explains Bishop's principle as
follows: ``Bishop criticized both non-constructive classical
mathematics and intuitionism.  He called non-constructive mathematics
`a scandal', particularly because of its `deficiency in {\bf numerical
meaning}'.  What he simply meant was that if you say something exists
you ought to be able to produce it, and if you say there is a function
which does something on the natural numbers then you ought to be able
to produce a machine which calculates it out at each number''
\cite{Fe00} [emphasis added--authors].}
The search for numerical meaning is a goal that can appeal to any
mainstream mathematician.  This thread will be pursued further in
Subsection~\ref{51}.

S.~Feferman \cite[end of section~A of part~I]{Fe00} identifies LPO as
a special case of the law of excluded middle, see Section~\ref{two}
above.

Given that a typical construction of Robinson's infinitesimals (see
Keisler \cite[p.~911]{Ke}) certainly does rely on LPO (but see
\cite{Pa}), Bishop's opposition to such infinitesimals, expressed in a
vitriolic%
\footnote{\label{vit}Historians of mathematics have noted the
vitriolic nature of Bishop's remarks, see e.g.,
Dauben~\cite[p.~139]{Da96}, cf.~footnotes~\ref{gab}, \ref{radical},
and \ref{gab2}.  M.~Artigue \cite[p.~172]{Ar} described it as
``virulent''; Davis and Hauser \cite{Dav78}, as ``hostile''; D.~Tall
\cite{Tal01}, as ``extreme''.}
review of Keisler's textbook, may have been expected.  Indeed, Bishop
wrote:
\begin{quote}
The technical complications introduced by Keisler's approach are of
minor importance.  The real damage lies in [Keisler's] {\bf
obfuscation\/} and devitalization of those wonderful ideas [of
standard calculus].  No invocation of Newton and Leibniz is going to
justify developing calculus using axioms V* and VI*-on the grounds
that the usual definition of a limit is too complicated!%
\footnote{Bishop is referring to the extension principle and the
transfer principle of non-standard analysis.  See
Subsection~\ref{hypo}.}
\cite[p.~208]{Bi77} [emphasis added--authors].
\end{quote}
M.~Davis notes that Bishop fails to acknowledge explicitly in his
review in the {\em Bulletin\/} that his criticism is motivated by his
foundational preoccupation with the law of excluded middle.%
\footnote{\label{f4}M. Davis writes as follows: ``Keisler's book is an
attempt to bring back the intuitively suggestive Leibnizian methods
(but see footnote~\ref{f49}) that dominated the teaching of calculus
until comparatively recently, and which have never been discarded in
parts of applied mathematics. A reader of Errett Bishop's review of
Keisler's book would hardly imagine that this is what Keisler was
trying to do, since the review discusses neither Keisler's objectives
nor the extent to which his book realizes them. Bishop[, meanwhile,]
objects to Keisler's description of the real numbers as a convenient
fiction (without informing his readers of the {\bf constructivist
context} in which this objection is presumably to be understood)''
\cite[p.~1008]{Dav} [emphasis added--authors].  See further
Subsection~\ref{hypo}.}

A similar point is alluded to in Keisler's brief (and measured)
response:
\begin{quote}
why did Paul Halmos, the {\em Bulletin\/} book review editor, choose a
constructivist as the reviewer?  \cite[p.~269]{Ke77}
\end{quote}
Keisler traced the source of Bishop's criticism, to the constructivist
criticism of classical mathematics, more specifically of its reliance
on classical logic, incorporating the key item rejected by
intuitionists, namely the law of excluded middle.  Comparing the use
of the latter to wine, Keisler compares Halmos' choice of reviewer, to
``choosing a teetotaller to sample wine''.  Halmos' reply will be
analyzed in Subsection~\ref{hal}.

In short, the constructivist Bishop is criticizing apples for not
being oranges: the critic and the criticized are not operating in a
common foundational framework.%
\footnote{A similar point was mentioned by M.~Davis~\cite{Dav}, see
footnote~\ref{f4}.}
The foundational framework of non-standard analysis, namely the
Zermelo-Fraenkel%
\footnote{\label{zf2}It is interesting to note a criterion of success
of a theory of infinitesimals as proposed by Adolf Abraham Fraenkel
and, before him, by Felix Klein.  In 1908, Klein formulated a
criterion of what it would take for a theory of infinitesimals to be
successful.  Namely, one must be able to prove a mean value theorem
for arbitrary intervals, including infinitesimal ones
\cite[p.~219]{Kl08}.  In 1928, A.~Fraenkel \cite[pp.~116-117]{Fran}
formulated a similar requirement in terms of the mean value theorem.
Such a Klein-Fraenkel criterion is satisfied by the
Hewitt-\Los-Robinson theory by the transfer principle, see
Appendix~\ref{rival}.}
set theory with the axiom of choice (ZFC), is the framework of the
vast majority of the readers of the {\em Bulletin}, at variance with
Bishop's preferred intuitionistic framework.  From the vantage point
of the latter, Bishop's {\em debasement of meaning\/} and {\em
obfuscation\/} charges would apply equally well to most of classical
mathematics, a point alluded to in Feferman's comment cited above.%
\footnote{See footnote~\ref{ff4} for Feferman's comment.}
Perhaps one may be permitted to detect a deficiency in the attributes
mentioned in the title of our current Section~\ref{kinder}, when
analyzing Bishop's published remarks concerning Robinson's approach.

\subsection{Halmos's translation}
\label{hal}

An interesting exchange took place in the pages of the {\em
Bulletin\/} and the {\em Notices\/}.  Keisler's response to Bishop's
review of his book was discussed in Subsection~\ref{eb}.  Halmos'
reply%
\footnote{\label{gab}An additional reply to Keisler's query, by
G. Stolzenberg \cite{Sto78b} (a close associate of Bishop's), occupies
slightly over a column on page 242 of the Notices.  Given its author's
interest \cite{Sto01} in literary deconstruction, we note that the
reply manages to employ a root, {\em dogma\/} (conspicuously absent
from Keisler's letter), on five occasions, culminating in the
expression {\em the spouting of dogma\/} in the penultimate paragraph,
cf.~footnotes~\ref{vit} and \ref{radical}.}
to Keisler's question came in the form of an editorial
pointer on p.~271 of the same issue, referring the reader to Halmos'
outline of his editorial philosophy on p.~283:
\begin{quote}
As for judgments, the reviewer may [\dots] say (or imply) what he
thinks.
\end{quote}
Halmos's philosophy was that a reviewer may legitimately use the
review of a book as a springboard for developing his own ideological
agenda.%

According to a close associate \cite{Ew} of Halmos', Halmos' policy
was to confront opposing philosophies in the goal of livening up the
debate.  One of his goals was to boost lagging sales that were
plaguing the publisher at the time, see \cite{Ha}.  The bottom-line
issue, combined with Halmos' own unflattering opinion of non-standard
analysis as ``a special tool, too special'' \cite[p.~204]{Ha},
apparently made the choice of Halmos' student (Bishop) as the
reviewer, attractive to the editor.

In the case of Keisler's book, such editorial philosophy translated
into a review, by Bishop, whose thinly disguised foundational agenda
(an attempt to de-legitimize the use of the law of excluded middle in
mathematical practice) took the form of vitriolic criticism of
Robinson-style infinitesimal calculus.

In his autobiography, Halmos described non-standard analysis as {\em a
special tool, too special\/} \cite[p.~204]{Ha}.  In fact, his
anxiousness to evaluate Robinson's theory may have involved a conflict
of interests.  In the early 1960s, Bernstein and Robinson~\cite{BR}
developed a non-standard proof of an important case of the invariant
subspace conjecture of Halmos'.  In a race against time, Halmos
produced a ``standard translation'' of the Bernstein-Robinson
argument, in time for the translation \cite{Ha2} to appear in the same
issue of {\em Pacific Journal of Mathematics\/}, alongside the
original.  Halmos invested considerable emotional energy (and
``sweat'', as he memorably puts it in his autobiography \cite{Ha})
into his translation.  Whether or not he was capable of subsequently
maintaining enough of a detached distance in order to formulate an
unbiased evaluation of non-standard analysis, his blunt unflattering
comments may create an impression of impropriety, as if he were
retroactively justifying his translationist attempt to deflect the
impact of one of the first spectacular applications of Robinson's
theory.

\subsection{Why was non-standard analysis the target?}
\label{whyfive}

We saw that Bishop's {\em debasement\/} charge applies to {\em all\/}
of classical mathematics.%
\footnote{\label{all1}See also footnotes~ \ref{all2a}, \ref{all2},
and~\ref{all3}.}
What was it specifically about non-standard analysis, among other
fields of classical mathematics, that may have prompted E.~Bishop to
speak out against it specifically?  We argue in Subsection~\ref{hypo}
that Bishop's comments on non-standard analysis in \cite{Bi75} and
\cite{Bi77} need to be understood in the context of his foundational
writings that appeared elsewhere, as neither his motivation, not his
terminology, are explained in \cite{Bi75} and \cite{Bi77}.  

We will first mention three possible technical reasons for Bishop's
choice of target, and then discuss a possible philosophical reason.
The possible technical reasons as the following:
\begin{enumerate}
\item
Certain results of classical analysis can be convincingly adjusted to
increase their numerical content, by strengthening their hypotheses
and/or weakening the conclusions (cf.~Subsection~\ref{shesh});
meanwhile, Bishop apparently felt that such an option is unavailable
for Robinson-style infinitesimals, leading to his wholesale rejection
of the latter.
\item
A related point is that the hyperreal approach incorporates an element
of non-constructivity at the basic level of the very number system
itself.
\item
By replacing Weierstrassian epsilontics%
\footnote{\label{f20}A term used by Dauben \cite{Da96}.  The
description of ``epsilon, delta'' arguments as {\em epsilontics\/}
(alternatively, epsilonics) is frequently found in the mathematical
literature, and even in the names of courses of instruction at certain
North American universities.}
by a theory of infinitesimals, non-standard analysis thereby removes
the bread and butter of constructive analysis, starting with the
constructive epsilontic definition of a Cauchy sequence, which
involves an explicit rate of convergence, see \cite{Bi67}.
\end{enumerate}

Recall that Brouwer sought to incorporate a theory of the continuum as
part of intuitionistic mathematics, by dint of his free choice
sequences.  Bishop's commitment to integr-ity (see glossary in
Subsection~\ref{glossary}) led to his dismissal of Brouwer's continuum
as a ``semimystical theory''.%
\footnote{\label{thicket}See main text at footnote~\ref{bugaboo} for
the full quote.}
Such a commitment on Bishop's part could hardly have countenanced an
about-face in the matter of the far denser thicket of Robinson's
hyperreal continuum.

More fundamentally, we would like to argue that non-standard analysis
presents a formidable philosophical challenge to constructivism, which
may, in fact, have been anticipated by Bishop himself in his
foundational speculations, as we explain below.

Bishop's vision of mathematics is at variance both with Kronecker's
constructivist program (see Subsection~\ref{kron}) and with Brouwer's
intuitionism.  Thus, Bishop's constructive mathematics claims to be
uniquely concerned with finite operations on the integers.  Meanwhile,
Bishop himself speculated that ``the primacy of the integers is not
absolute'':
\begin{quote}
It is an {\bf empirical fact\/} that all [finitely performable
abstract calculations] reduce to operations with the integers.  There
is no reason mathematics should not concern itself with finitely
performable abstract operations of other kinds, in the event that such
are ever discovered \cite[p.~53]{Bi68} [emphasis added--authors].
\end{quote}
Bishop hereby acknowledges that the {\em primacy of the integers\/} is
merely what he describes as an {\em empirical fact\/}.  The
implication is that the said {\em primacy\/} could be contradicted by
novel mathematical developments.  Perhaps Bishop sensed that a
rigorous theory of infinitesimals is both
\begin{itemize}
\item
not reducible to finite calculations on the integers, and yet
\item
accomodates a finite performance of abstract operations,
\end{itemize}
thereby satisfying his requirements for legitimate mathematics.
Bishop made a foundational commitment to the {\em primacy of the
integers\/} (a state of mind known as {\em integr-ity\/}, see
Subsection~\ref{glossary}) through his own work and that of his
disciples starting in the late sixties.  He may therefore have found
it quite impossible, in the seventies, to acknowledge the existence of
``finitely performable abstract operations of other kinds''.  We will
analyze the related issue of post-LEM numerical meaning in
Subsection~\ref{postlem}.

\subsection{Galley proof ``debasement''}
\label{debasement}

It has recently come to light that Bishop's ``debasement'' remark,
published in his {\em Crisis\/} text \cite{Bi75}, in reality was never
uttered in his oral presentation in the presence of the distinguished
audience%
\footnote{See footnote~\ref{distinguish}.}
at the 1974 workshop in Boston.  A careful reading of Bishop's {\em
Crisis\/} text already suggests that his ``debasement'' comment may
have been added at the galley proof stage, as its deletion from the
text causes no interruption of textual continuity.  In fact, the
second-named author has been able to locate a participant in the
workshop, who confirmed that Bishop said nothing about non-standard
analysis in his oral presentation.  The participant in question, a
historian of mathematics, wrote as follows on the subject of Bishop's
statement published in 1975 on Robinson-style infinitesimal calculus:
\begin{quote}
I do not remember that any such statement was made at the workshop and
doubt seriously that it was in fact made.  I would have pursued the
issue vigorously, since I had a particular point of view about the
introduction of non-standard analysis into calculus.  I had been
considering that question somewhat in my attempts to understand
various standards of rigor in mathematics.  The statement would have
fired me up \cite{Man}.
\end{quote}
Furthermore, the text~\cite{KS}, entitled {\em Commentary on
E.~Bishop's talk\/}, states clearly that the said commentary was based
on a {\em transcript\/} of Bishop's actual talk at the workshop.
%
%
When the second-named author contacted the authors of the said {\em
Commentary\/} in the summer of~2009, it turned out that neither was
aware of the fact that the {\em published\/} version of Bishop's {\em
Crisis\/} essay contained a paragraph on the subject of non-standard
analysis.%
\footnote{\label{gab2}In fact, Stolzenberg replied on 23 june 2009 to
the effect that he is ``morally certain [that Bishop] didn't mention
non-standard analysis during that lecture''.  The second-named author
subsequently requested a copy of the {\em transcript\/} of Bishop's
lecture (on which the {\em Commentary\/} \cite{KS} was based).  The
reply arrived on 17 august 2009, and consisted of three parts:
\begin{enumerate}
\item
it acknowledged that the authors are in possession of such a
transcript;
\item
it denied the request for a copy of the transcipt;
\item
it expressed concern about ``the harm'' that would result if this
writer were to obtain a copy of the trascript.
\end{enumerate}
The tenor of these remarks suggests that in fact the transcript of
Bishop's talk contains no mention of non-standard analysis whatsoever,
which would therefore prove that the addition of the vitriolic
paragraph was an afterthought.  The community of historians of
mathematics should encourage Stolzenberg to release the transcript of
Bishop's lecture in Stolzenberg's possession, as it represents an
important historical document.}
Thus, Bishop appears to have added his ``debasement'' comment on
non-standard calculus at the galley proof stage of publication.

To be sure, additional material is frequently added at galley stage.
Was the addition of the ``debasement'' comment ethically problematic?
Note the published version of Bishop's lecture includes an epilogue.
The epilogue contains the discussion (following the lecture),
involving Birkhoff, Mackey, J.P.Kahane, and others, who have
challenged Bishop on a number of points.  None of the participants
challenged Bishop on his criticism of infinitesimals \`a la Robinson.
This reader at least of Bishop's essay was disappointed, upon his
first reading of the essay, not to find any challenge to Bishop's
``debasement'' of non-standard analysis, in the discussion that
ensued.  We now understand why there was no such challenge: Bishop did
not say a word about Robinson in his oral presentation.

The important ethical point concerns the juxtaposition of Bishop's
vitriolic comment (added at galley stage) concerning non-standard
analysis, and an absence of any reaction to it in the ensuing
discussion.  Such a juxtaposition is misleading, as it suggests, if
not agreement, then at least a tolerance toward Bishop's comment on
the part of the participants in the discussion panel, raising issues
of {\em integrity\/} (without the dash).%
\footnote{Bishop should have been aware of the fact that the follow-up
discussion will be included in the published version of his lecture.
In order to avoid misunderstandings, he could have easily included a
footnote to the effect that the non-standard paragraph was added at
the galley proof stage, but he didn't.}
The vitriolic tenor of Bishop's remark was aimed ostensibly at
Robinson-style infinitesimal calculus, but logically would apply
equally well to the rest of mainstream mathematics, as we discuss in
the next subsection.%
\footnote{\label{all2a}See Section~\ref{hypo} and Feferman's comment
cited in footnote~\ref{ff4}, as well as footnotes~\ref{all1},
\ref{all2}, and~\ref{all3}.}

\subsection{A hypothesis and its refutation}
\label{hypo}

Bishop diagnosed classical mathematics with a case of a ``debasement
of meaning'' in his {\em Schizophrenia in contemporary mathematics\/},
written in 1973.%
\footnote{See Bishop \cite[p.~1]{Bi85}.}
Hot on the heels of {\em Schizophrenia\/} came the 1975 {\em Crisis in
contemporary mathematics\/}~\cite{Bi75}, where the same ``debasement
of meaning'' diagnosis was slapped upon infinitesimal calculus \`a la
Robinson.  It is therefore transparent that his opposition to
Robinson-style infinitesimals is merely an instance of a broad
opposition to classical mathematics as a whole, on
ideological/foundational grounds.

Nonetheless, a hypothesis has been advanced by a number of
constructivists%
\footnote{Including Stolzenberg \cite{Sto78b}.}
that Bishop's criticism of infinitesimal calculus \`a la Robinson was
motivated primarily by pedagogical, rather than foundational,
concerns.  According to such a hypothesis, his objection focused on
the use of the extension principle and the transfer principle of
non-standard analysis (see Appendix~\ref{rival}), whose acceptance by
the student would necessarily be on the grounds of trust in the
teacher.%
\footnote{\label{potential}It has been argued further that the idea of
constructive mathematics is to found mathematics on actions that
everyone can concretely perform, in such a way that mathematics is
accessible for everyone.  Let us examine the coherence of such an
argument (certain implications of a Bishop-style constructive approach
in teaching is examined in Subsection~\ref{shesh}).  The assumption of
a completed (actual) infinity of~$\Z$ in Bishopian constructivism is a
departure from Kronecker's vision.  It cannot be ``performed'' in any
concrete sense by either student or teacher.  Nevertheless, both
constructivists and classical mathematicians acknowledge that such an
assumption facilitates both teaching and research in mathematics.}

However, such a hypothesis does not stand up to scrutiny, as we argue
below.  A widespread perception of a foundational motivation behind
Bishop's criticism of non-standard analysis has been attested to by
M. Davis, V. Komkov, J. Dauben, H. J. Keisler, D. Tall, and others,
and is not a novel interpretation.

Thus, M.~Davis comments that Bishop states his objections ``without
informing his readers of the constructivist context in which this
objection is presumably to be understood'' \cite[p.~1008]{Dav}.
Physicist V.~Komkov notes that
\begin{quote}
Bishop is one of the foremost researchers favoring the constructive
approach to mathematical analysis.  It is hard for a constructivist to
be sympathetic to theories replacing the real numbers by hyperreals
\cite[p.~270]{Kom}.%
\footnote{\label{komk}Leaving aside the issue whether or not
non-standard analysis can be done constructively, Komkov's perception
of a foundational concern on Bishop's part is unmistakable.}
\end{quote}
Philosopher of mathematics G.~Hellman writes:
\begin{quote}
Some of Bishop's remarks (1967) suggest that his position belongs in
[the radical constructivist] category \cite[p. 222]{Hel93a}.
\end{quote}
Historian of mathematics J.~Dauben analyzed Bishop's criticism and
found it to be ``unfounded''.%
\footnote{\label{daub}After evoking the success of nonstandard
analysis ``at the most elementary level at which it could be
introduced--namely, at which calculus is taught for the first time'',
Dauben proceeds to point out that ``there is also a deeper level of
meaning at which nonstandard analysis operates.''  Dauben mentions the
impressive applications in ``physics, especially quantum theory and
thermodynamics, and in economics, where study of exchange economies
has been particularly amenable to nonstandard interpretation'' (some
references for such applications can be found in
Appendix~\ref{rival}).  At this deeper level of meaning, Dauben
concludes: ``Bishop's views can be questioned and shown to be as
unfounded as his objections to nonstandard analysis pedagogically''
\cite[p.~192]{Da88}.}

Through a comparative textual study of Bishop's three texts, we show
that his criticism of non-standard analysis was motivated more by
foundational than by pedagogical issues.

\medskip\noindent 1. Bishop's perception of a {\em shared\/}
shortcoming in non-standard calculus, on the one hand, and classical
mathematics pedagogy, on the other, is evident from his contention
that the former will
\begin{quote}
{\bf confirm\/} the students' previous experience \cite[p.~208]{Bi77}
[emphasis added--authors].
\end{quote}
of the latter.  The usage of the term {\em confirm\/} points to a {\em
shared\/} shortcoming.

\medskip\noindent
2. The nature of such a shared shortcoming is hinted upon in his
reference to the students' experience with classical mathematics as a
``{\bf meaningless\/} exercise in technique'' \cite[p.~208]{Bi77}
[emphasis added--authors].

\medskip\noindent 3. Note that the term {\em meaning\/} is alluded to
elsewhere in his review, as well.  Thus, he laments that in ``our
educational system'', it is ``a bad form to ask what it all means''
\cite[p.~206]{Bi77}.  The 1977 review is evasive as to what exactly
Bishop has in mind when he mentions {\em meaning\/}.

\medskip\noindent 4. Meanwhile, Bishop's {\em Crisis\/} essay
\cite{Bi75} from 1975 focuses on the virtue of {\em numerical
meaning\/}, defined ultimately as the rejection of the law of excluded
middle even in its special case called LPO.%
\footnote{See footnote~\ref{ff4} for a more detailed discussion of
Bishop's LPO.}
Bishop criticized non-standard calculus as a ``debasement of meaning''
in the following terms:
\begin{quote}
It is difficult to believe that debasement of meaning could be {\bf
carried so far} \cite[514]{Bi75} [emphasis added--authors].
\end{quote}
Bishop's wording suggests that {\em debasement of meaning\/} had
already been spoken of, while, allegedly, this particular instance of
it, is {\em carried\/} particularly {\em far\/}.  However, the {\em
Crisis\/} text does not contain any other occurrences of the
three-word formation.  Where was {\em debasement of meaning\/} spoken
of previously by Bishop?

\medskip\noindent 5. The realm of applicability of the three-word
formation is clarified in Bishop's essay {\em Schizophrenia in
contemporary mathematics\/} \cite{Bi85}.  The essay was published
posthumously in 1985, but M. Rosenblatt \cite[p.~ix]{Ros} points out
that it was ``distributed'' in 1973.  In other words, Bishop's {\em
Schizophrenia\/} text was distributed a year or two prior to the {\em
Crisis\/} text.  Bishop writes as follows:
\begin{quote}
Brouwer's criticisms of classical mathematics were concerned with what
I shall refer to as ``the debasement of meaning'' \cite[p.~1]{Bi85}.%
\footnote{\label{all2}See also footnotes~\ref{all1}, \ref{all2a},
and~\ref{all3}.}
\end{quote}
In Bishop's own words, the {\em debasement of meaning\/} expression,
employed in his {\em Crisis\/} text to refer to non-standard calculus,
is in fact a criticism of {\em classical mathematics\/} as a whole.%
\footnote{When Bishop's usage of the identical expression, {\em
debasement of meaning\/}, in the context of both non-standard calculus
and classical mathematics, was presented in 2009 (as evidence of the
claim of a foundational motivation on Bishop's part in criticizing
non-standard calculus) to a constructivist and a close associate of
Bishop's, he responded as follows: ``Classical math here; non math
major calculus there.  Two different kinds of debasement of meaning''.
See footnote~\ref{radical}.}
The three-word formation sailed effortlessly from {\em
Schizophrenia\/} to {\em Crisis\/}, revealing its author's intention.

\medskip\noindent
6.  By the time Bishop got around to describing non-standard calculus
as a {\em debasement of meaning\/} \cite[p.~514]{Bi75}, he had already
used the term {\em meaning\/} as shorthand for {\em numerical
  meaning\/} in his invocation of a philosophical principle, to the
effect that a discussion of {\em meaning\/} should precede a
discussion of {\em truth\/} \cite[p.~509]{Bi75}.

\medskip\noindent
7. Bishop's remark about meaning preceding truth occurs in the context
of his analysis of the Brouwer-Hilbert controversy, which he
attributes to the fact that
\begin{quote}
``[t]hey attached different {\bf meanings} to Cantorian set theory''
\cite[p.~509]{Bi75} [emphasis added--authors]
\end{quote}

\medskip\noindent
8. Bishop provides an example explaining the meaning of an integer to
Brouwer:
\begin{quote}
``[f]or example, an integer to Brouwer [\dots] is either a[ concrete]
integer in decimal notation or a method that in principle will lead
after a {\em finite\/} number of steps to an integer in decimal
notation.  Again there is this notion of computability: if the integer
is not given directly to be sure that it is finitely {\bf
computable}.''  \cite[p.~509-510]{Bi75} [emphasis added--authors]
\end{quote}

\medskip\noindent
9. Toward the end of an imaginary dialog between Brouwer and Hilbert,
Bishop has Brouwer express himself as follows on the subject of
meaning:
\begin{quote}
For me to do my mathematics in your system would entail a significant
loss of meaning%
\footnote{Bishop does not speculate what Hilbert's reply might have
been to such a charge of a loss of meaning, but we cannot help noting
that for a classical mathematician to do mathematics in Brouwer's
system, would entail a significant loss of post-LEM numerical meaning.
Thus, H.~Weyl, who became disillusioned with the Brouwerian
development by the 1930s, notes that ``Mathematics with Brouwer gains
its highest intuitive clarity. He succeeds in developing the
beginnings of analysis in a natural manner, all the time preserving
the contact with intuition much more closely than had been done
before. It cannot be denied, however, that in advancing to higher and
more general theories the inapplicability of the simple laws of
classical logic eventually results in an almost unbearable
awkwardness. And the mathematician watches with pain the larger part
of his towering edifice which he believed to be built of concrete
blocks dissolve into mist before his eyes'' \cite[p. 54]{We49} (see
further in Feferman \cite{Fe97}).  See also Subsection~\ref{postlem},
and particularly Example~\ref{action}.}
\cite[p.~511]{Bi75}.
\end{quote}
To Brouwer, as interpreted by Bishop, {\em meaning\/} is seen through
the lens of notions of {\em computability\/}, {\em numerical
meaning\/}, and LEM.  Meanwhile, Hilbert is uninterested in LEM
extirpation; for him, the source of meaning lies elsewhere.%
\footnote{\label{hilbert}Thus, Novikov \cite{No2} writes that, for all
of Hilbert's interest in a program of formalisation and
axiomatisation, he understood such a program in a highly nontrivial
fashion.  Novikov mentions, notably, Hilbert's deep theorem
\cite{Hi15} to the effect that the Einstein equations of relativistic
gravitation are of Lagrangian type.  Novikov's essay is discussed
further in Section~\ref{fifteen}.  A detailed discussion of the
significance, and {\em meaning\/}, of Hilbert's program may be found
in Avigad and Reck \cite{AR}.  See also Corry's quote of Hilbert
(footnote~\ref{corry}).}

\medskip\noindent 10. Thus, the shared shortcoming of infinitesimal
calculus \`a la Robinson and classical mathematics, in Bishop's eyes,
ultimately boils down to a deficiency in numerical meaning, in other
words, to a reliance on LEM.

\medskip
Through a comparative textual analysis of his three texts, we have
argued that Bishop's grievance is foundationally motivated, even when
couched in a pedagogical idiom.  Indeed, foundational ideology and
pedagogy were inseparable in Bishop's thinking.

Our conclusion is that the Bishop-Keisler controversy was, in essence,
a foundational controversy, as suggested by J. Dauben, see
\cite[p.~177]{Da88}, \cite[p.~190]{Da88}, and \cite[p.~192]{Da88},%
\footnote{See footnote~\ref{daub}.}
as well as by V.~Komkov.%
\footnote{See footnote~\ref{komk}.}
The controversy involved a fundamental disagreement concerning a basic
question, which can be formulated in terms of the dichotomy introduced
in Subsection~\ref{postlem}, as follows.  Is {\em post-LEM numerical
meaning\/}, mathematically coherent?

In conclusion, Bishop's vitriolic criticism of infinitesimals \`a la
Robinson was but a salvo in his broadside attack on classical analysis
and classical mathematics as a whole, percolating through a number of
his foundational texts, and culminating in the posthumously published
charge of {\em Schizophrenia\/} \cite{Bi85}.

\subsection{Further background on refuted hypothesis}
\label{hypo2}

In the previous section, we argued that Bishop's criticism
\cite{Bi75}, \cite{Bi77} of non-standard calculus was foundationally
motivated.  Here we provide some background on the context of Bishop's
review, and clarify some of the pedagogical and mathematical issues
involved.

Bishop's {\em Crisis\/} \cite{Bi75} tells an epic tale of how, half a
century earlier, a failure of communication between the two giants
(Brouwer and Hilbert), had, in Bishop's view, the effect of derailing
mathematics from its numerical foundations, and set it on its present
``meaningless'' course.  {\em Crisis\/} also presents a vision of how
mathematics may have been rescued, had the Brouwer-Hilbert debate only
followed a different script (namely, the one specified
in~\cite{Bi75}), focusing on the foundational issue of ``principles of
omniscience'' and the problem of its acceptability.

In the midst of the epic tale, Bishop inserts%
\footnote{\label{f22}See Subsection~\ref{debasement}.}
a vitriolic paragraph against Robinson-style infinitesimal calculus.
To believe Bishop's followers, his unique concern is the best way to
teach freshman calculus.  While a certain suspension of disbelief is
required here, it should also be pointed out that if the source of
Bishop's interest were indeed pedagogical, then he pursued it with
extraordinary nonchalance.

The remarks he added to the galley proofs$^{\ref{f22}}$ of \cite{Bi75}
indicate that he was already aware of infinitesimals in the classroom,
at least two years prior to writing his 1977 book review~\cite{Bi77}.
Meanwhile, K.~Sullivan's 1976 field study \cite{Su} of infinitesimals
in the classroom was published in a mainstream periodical, was easily
accessible, and showed the infinitesimal approach to be more
successful than the standard one in imparting an ability to interpret
the sense of mathematical formalism of calculus.  J. Dauben pointed
out that
\begin{quote}
[Sullivan's] study [\dots] was presumably available to Bishop when his
review of Keisler's book appeared in 1977, in which he attacked the
pedagogical validity of nonstandard analysis \cite[p.~190]{Da88}.
\end{quote}
Dauben analyzes Sullivan's comparison of the nonstandard approach and
the traditional approach, with the conclusion that
\begin{quote}
contrary to Bishop's views, the traditional [epsilontic] approach to
calculus may be the more pernicious \cite[p.~191]{Da88}.
\end{quote}
Why wasn't Sullivan's study dealt with, one way or another, in
Bishop's book review?  Could it be because Bishop was not interested
in field studies?  Could his opposition to the non-standard approach
have been based, not on any of its perceived shortcomings in imparting
knowledge of rates of change, derivatives, areas, integrals, etc., but
rather on ideological and foundational grounds?  That a certain amount
of foundational material should be expected to be taken on trust in a
first calculus course, is a feature that infinitesimal-based calculus
shares with the standard approach.  Very rare indeed is the first
calculus course that deals with equivalence classes of Cauchy
sequences (or other constructions of the real number field).  In the
decade of his life that Bishop seems to have become interested in
calculus methodologies, he is not known to have written critical book
reviews of any methodologies other than the hyperreal one,
which--curious coincidence!--incorporates an element of
non-constructivity at the basic level of the very number system
itself.%
\footnote{See footnote~\ref{f7} for Robinson's take on number
systems.}

Having already written, about Bertus Brouwer himself, that
\begin{quote}
``Brouwer's bugaboo has been compulsive speculation about the nature
of the continuum.  His fear seems to have been that, unless he
personally intervened to prevent it, the continuum would turn out to
be discrete.  [The result was Brouwer's] semimystical theory of the
continuum''%
\footnote{\label{bugaboo}See also main text at
footnote~\ref{thicket}.}
\cite[p.~6 and 10]{Bi67},
\end{quote}
was Bishop likely to take lying down, a far denser thicket of
Robinson's hyperreals?  Could his interest have been foundationally
motivated, as indeed it was {\em perceived\/} by many of the readers
of the {\em Bulletin\/}, as is readily acknowledged even by Bishop's
followers?

To elaborate, note that Bishop quotes Keisler's observation that 

\begin{quote}
``we have no way of knowing what a line in physical space is really
like.  It might be like the hyperreal line, the real line, or neither.
However, in applications of the calculus, it is helpful to imagine a
line in physical space as a hyperreal line.''  
\end{quote}

Bishop himself does not comment, other than voicing a criticism as to
an alleged lack of supporting evidence.  His presentation suggests
that Bishop himself favors the classical real line.  Yet Bishop states
clearly in his {\em Manifesto\/} \cite{Bi67} that from his
constructivist viewpoint, it is LEM, and not the axiom of choice, that
is the source of non-constructivity.%
\footnote{See our main text around footnote~\ref{lessthan}.}

To summarize, Bishop seeks to convince his classical reader of the
{\em idealistic\/} nature (see glossary in Section~\ref{glossary}) of
a Robinson infinitesimal, whereas to Bishop himself, the classical
continuum%
\footnote{He felt similarly about Brouwer's continuum; see
footnote~\ref{kronecker} and the main text at footnote \ref{bugaboo}.}
should not ``have been put there in the first place'' any more than a
Robinson infinitesimal, as the classical construction of both relies
on the law of excluded middle.  The {\em similar\/} foundational
status of the classical continuum, on the one hand, and the hyperreal
number, on the other, constitutes an astonishing point of convergence
between Bishop and Robinson.%
\footnote{\label{f7}Concerning the reality of number systems, Robinson
expressed himself as follows: ``[T]he infinitely small and infinitely
large numbers of a non-standard model of Analysis are neither more nor
less real than, for example, the standard irrational numbers [\dots]
both standard irrational numbers and non-standard numbers are
introduced by certain infinitary processes" \cite[p.~282]{Ro66}.}

On the subject of Keisler's infinitesimals, Bishop had this to say:
\begin{quote}
We are not told what an infinitesimal~$\Delta x$ is [\dots] Keisler
has developed limits from a supposedly consistent system of axioms
[\dots] But he has not explained the axioms.  They are mere
conveniences for generating proofs, whose intuitive content will {\bf
certainly} escape the students \cite[p.~207]{Bi77} [emphasis
added--authors].
\end{quote}

What exactly is the source of Bishop's {\em certainty\/}?  The point
is not so much that, having taught calculus effectively, first-hand,
from Keisler's book, this author disagrees with Bishop's assessment as
to the effectiveness of Keisler's method.%
\footnote{See \cite{KK10a, KK10b} for details.}
Rather, the point is that the nature, and the language, of Bishop's
criticisms clearly reveal an axiomatic/foundational preoccupation.

Dauben~\cite{Da96} also understood Bishop's ``debasement of meaning''
criticism of non-standard calculus as referring to an alleged lack of
numerical meaning.%
\footnote{See footnote~\ref{daub}.}
The latter is Bishop's term for the {\em sine qua non\/} ingredient of
constructively acceptable mathematics.  The lack of numerical meaning,
to Bishop, corresponds to a reliance on LPO or more generally, the law
of excluded middle (see glossary in Subsection~\ref{glossary}).

The over-reliance on the teacher's authority, in Bishop's view, is
only a shortcoming if such authority is propped up by a
non-constructive arsenal such as excluded middle or classical choice.
Bishop expresses the sentiment that
\begin{quote}
the notorious~$\epsilon, \delta$ definition of limit is {\bf common
sense} \cite[p.~208]{Bi77} [emphasis added--authors].
\end{quote}
He further shares with the reader the following confidence: ``[the
students] don't believe me''.  Dauben speculates that there may be a
very good reason for their disbelief, namely,
\begin{quote}
that what [Bishop] claims, in fact, does not seem to be true
\cite[p.~134]{Da96}.
\end{quote}
The vast education literature on the student reception of quantifiers
(see a useful bibliography in \cite{WAD}) attests to the irreducible
difficulty of properly inculcating an understanding of even a single
alternation of quantifiers, let alone an iteration of such
alternations involved in the~$\epsilon, \delta$ definition.%
\footnote{A detailed discussion of a typical case study appears in the
main text at footnote~\ref{fernandez}.}

However, our main goal here is not to gauge the difficulty of
Weierstrassian definition of the limit concept based on epsilontics.
Rather, the point is to signal the nonchalance of Bishop's comments,
if they are to be interpreted as pedagogically motivated.

Interpreted as foundationally motivated, on the other hand, Bishop's
comments become intelligible.  In traditional (LEM-dominated)
mathematics courses, students become accustomed to going through the
numerically meaningless motions of {\em idealistic mathematics\/} (see
glossary in Section~\ref{glossary}), so that by the time they reach
Bishop's classroom and are exposed to the numerically meaningful
content of Weierstrassian epsilontics, they are put off by its
pristine ``common sense''.  As part of a constructivist critique of
the foundations, Bishop's criticism is intelligible and compelling,
even if ultimately flawed, as we argue in Section~\ref{fourteen}.

For all its touted advantages of increasing numerical meaning by means
of a constructivisation of all that can be constructivized, Bishop's
approach is not without its pitfalls.  Note that a ``de-{\em
basement\/}" of mainstream mathematics is what constructivism is
sometimes attempting to carry out, by changing the meaning of {\em
basic\/} terms such as {\em continuous function\/} and {\em finite
set\/}, as analyzed in Sections~\ref{shesh} and \ref{six}.

\section{De-basing in constructivism}
\label{debasing}

A number of familiar terms in the mathematical lexicon are displaced
from their base, or basic, meaning when they are rendered
constructive.  We will analyze several such phenomena in this section.

\subsection{Verbal equality}
\label{shesh}

Bishop's own approach to constructive mathematics only accepts
uniformly continuous functions (more precisely, uniformly continuous
on compact intervals); they are accordingly renamed {\em continuous\/}
in Bishop's foundational text~\cite{Bi67}.%
\footnote{While Bishop, at the time of writing \cite{Bi67}, appears to
have thought that such functions were the only ones that would enter
into numerically meaningful mathematics, subsequent research has shown
that even sequential continuity (constructively strictly weaker than
even pointwise continuity) arises naturally in constructive functional
analysis; see Ishihara \cite{Is}.}
S.~Feferman writes that Bishop deals with
\begin{quote}
a very special class of functions of real numbers, namely those which
are uniformly continuous on [e]very compact interval.  In this way, he
finesses the whole issue of how one arrives at Brouwer's theorem [to
the effect that every function on a closed interval is uniformly
continuous]%
\footnote{A discussion of the issue of the finessing of Brouwerian
counterexamples in the context of the extreme value theorem appears in
Subsection~\ref{evt}.}
by saying that those are the only functions, at least initially, that
one is going to talk about \cite{Fe00}.
\end{quote}

Continuity differs from uniform continuity in at least the following
three ways:

\begin{enumerate}
\item
continuity can be defined in terms of limits; 
\item
continuity is a pointwise concept;
\item
only~$3$ quantifiers are required to define continuity at a point.
\end{enumerate}

To elaborate, note that a typical calculus textbook \cite[p.~108]{TF}
defines continuity in terms of a previously defined notion of the
limit.  This is possible because ordinary continuity is a local
(pointwise) property, whereas to define uniform continuity, one needs
to work with pairs of points.  Thus, the definition of uniform
continuity \cite[p.~115]{TF} requires a return to~$\epsilon,
\delta$'s.

Most textbooks adopt the standard Weierstrassian approach using
epsilontics.  When it comes to continuity and uniform continuity, both
definitions involve alternations of quantifiers.  However,
pedagogically speaking, one can present the definition of ordinary
continuity with only three quantifiers, once the point has been fixed.
There is no such device that would permit a definition of uniform
continuity with fewer than four quantifiers.  Such a definition is
almost inaccessible for an average undergraduate.

A revealing case study \cite{Fer} documents an effort properly to
explain~$\epsilon, \delta$ to students at ``a state university with an
enrollment of 15,000'', using an interactive approach that reportedly
generated considerable interest on the part of the students.  It
emerges that due to the intrinsic difficulty of the subject matter,
the semester was spent explaining~$\epsilon, \delta$ for {\em
linear\/} functions.  The course itself never reached quadratic
functions, but a number of students ``began coming to office hours on
a regular basis'' to get a taste of the {\em quadratic\/} case.  The
author concludes that
\begin{quote}
[a]lthough the ideas concerning quadratics were pursued outside of
class with only a limited number of students, it is important to note
that the students initiated this discussion and pursued it on their
own time \cite[p.~53]{Fer}.%
\footnote{\label{fernandez}There are indeed numerous case studies of
this sort; see e.g. \cite{WAD}.}
\end{quote}
For educators, it is no mystery why Johnny can't take limits (see
Roh~\cite{Roh}).  Many teachers, based on classroom experience, feel
that uniform continuity is a concept that is even more difficult for
students to grasp than ordinary continuity.

Bishop's renaming of uniformly continuous functions as {\em
continuous\/}, is part of a pattern commented upon by reviewer B. van
Rootselaar (for Mathscinet) of Bishop's book~\cite{Bi67}, who writes
as follows:

\begin{quote}
The comparison with classical mathematics carried out by the author in
the notes accompanying the chapters is somewhat superficial.  In fact
the author stresses the equal hypothesis interpretations of classical
theorems, which is misleading since the hypotheses are for the greater
part only verbally equal \cite{vRoo}.
\end{quote}

What van Rootselaar is referring to is the fact that some of the
theorems are worded in such a way as to sound identical to the
classical ones, but, of course, the definitions have been tampered
with, e.g.~what is defined as {\em continuity\/} is really {\em
uniform continuity\/}, classically speaking (See also Hellman
\cite{He89}).

Apart from the issue of {\em verbal equality\/}, the numerical
advantage of Bishop's approach has been neatly summarized by Feferman,
who pointed out that in constructive analysis,
\begin{quote}
concepts are chosen so that there is a lot of witnessing information
introduced in a way that is not customary in classical mathematics,
where it is hidden, for instance, by implicit use of the axiom of
choice \cite{Fe00}.
\end{quote}

\subsection{Goodman-Myhill's 
$1\frac{1}{2}\pm \frac{1}{2}$ sets and ``full'' choice}
\label{six}

What does the term {\em finite\/} mean?  In this section we will
comment on the constructive distinction between finite and subfinite
sets, in the context of a result relating the axiom of choice and LEM.

A theorem dating from 1975, alternatively referred to as the the
Diaconescu theorem \cite{Di} or the Goodman-Myhill theorem \cite{GM},
asserts that the axiom of choice implies the law of excluded middle.

There is an intriguing point here that Bishop actually claimed that he
has no quarrel with the axiom of choice.  In his {\em Constructive
manifesto\/}, eight years prior to the publication of
Diaconescu-Goodman-Myhill, he wrote:
\begin{quote}
the axiom of choice [\dots] is not a real source of the
unconstructivities of classical mathematics [\dots] [it] is used to
extract elements from equivalence classes where they should never have
been put in the first place \cite[p.~9]{Bi67}.
\end{quote}
His remark has apparently been a source of discomfort for some of his
followers.%
\footnote{\label{lessthan}It appears to some, if not all, constructive
interpreters that Bishop's remark was, if studied carefully, perfectly
in order, if less than perfectly expressed; see \cite[p.~13]{BB}.  
%
%
}
Understanding the proof of the Diaconescu-Goodman-Myhill theorem
hinges on the following critical distinction between finite and
subfinite sets.%
\footnote{The constructivist position on the issue can be summarized
as follows.  The contention that the image of a finite set under a
function~$f$ need not be finite, has precise meaning in a constructive
context, reflecting a deficiency of numerical meaning in a classical
assertion of its finiteness.}
A set possessing either one or two elements (we don't know which)
cannot be called {\em finite\/} in the constructive sense.  Finding a
choice function for a collection of sets containing either one or two
sets (we don't know which) requires a nontrivial application of the
``full" axiom of choice.  The theorem is proved by means of such an
application of the axiom of choice (see Subsection~\ref{two} for a
discussion of the nature of the quantifiers in the axiom of choice).%
\footnote{For the benefit of a classically-trained reader it may be
helpful to clarify that the ``fullness'' of the power of the axiom of
choice referred to here is not in the usual sense, of applying it to a
collection of sets of arbitrarily high cardinality, but rather to the
{\em omniscience\/} aspect of the axiom of choice (see discussion of
LPO in Subsection~\ref{eb}), stemming from the classical
interpretation the the existence quantifier.}
Note that we are dealing with a collection of sets
containing~$1\frac{1}{2}\pm \frac{1}{2}$ sets (namely, either~$1$
or~$2$ sets).%
\footnote{The distinctions made by constructive mathematicians are
real and important if you want to work constructively.  Yet they often
look odd if you cannot step outside the classical framework.}

Note that the law of excluded middle is a constructive consequence of
much more elementary statements than that of the full axiom of
choice.%
\footnote{\label{bm} Let~$P$ be any statement, and let~$S=\{0\} \cup
\{x: x = 1 \wedge (P \vee \neg P)\}$.  Then~$S$ is a subset of
$\{0,1\}$ that certainly contains~$0$; but~$1 \in S$ if and only if
$P\vee \neg P$ holds.  Thus if {\em under the
Brouwer-Heyting-Kolmogorov interpretation of logic\/}, we could prove
that every inhabited subset of~$\{0,1\}$ either contains exactly one
element, or contains both~$0$ and~$1$, then we could prove the law of
excluded middle.  So the law of excluded middle is a constructive
consequence of much more elementary statements than that of the full
axiom of choice.
%
%
}

\subsection{Is Bishopian constructivism compatible with classical logic?}
\label{evt}

We will examine the question of such compatibility in this subsection,
in the context of the extreme value theorem.  Note that the extreme
value theorem does not hold constructively.  Namely, the existence of
a maximum of a continuous function~$f$ on~$[0,1]$, in the sense of the
{\em constructive\/} formula
\[
\left(\exists_{\text{InLo}}\, x \in [0,1]\right) \; \left( f(x)=\sup(f)
\right),
\]
does not hold.  Indeed, let~$a$ be any real number such that
\[
a\leq 0 \text{\;or\;} a\geq 0 \text{\;is unknown} .
\]
Next, define~$f$ on~$[0,1]$ by
\[
f(x):=ax.
\]
Then~$\sup(f)$ is simply~$\max(0, a)$, but the point where it is
attained cannot be captured constructively (see~\cite[p.~295]{TV}).
To elaborate on the foundational status of this example, note that the
law of excluded middle:~$P \vee \neg P$ (``either~$P$ or (not~$P$)''),
is the strongest principle rejected by constructivists.  A weaker
principle is the LPO (limited principle of omniscience).  The LPO is
the main target of Bishop's criticism in \cite{Bi75}.  The LPO is
formulated in terms of sequences, as the principle that it is possible
to search ``a sequence of integers to see whether they all vanish''
\cite[p.~511]{Bi75}.  The LPO is equivalent to the law of
trichotomy:
\[
(a< 0) \vee (a=0) \vee (a> 0).
\]
An even weaker principle is
\[
(a\leq 0) \vee (a \geq 0),
\]
whose failure is exploited in the construction of the counterexample
under discussion.  This property is false intuitionistically.  After
discussing real numbers~$x\geq 0$ such that it is ``{\em not\/}'' true
that~$x>0$ or~$x=0$, Bishop writes:
\begin{quotation}
In much the same way we can construct a real number~$x$ such that it
is {\em not\/} true that~$x\geq 0$ or~$x\leq 0$ \cite[p.~26]{Bi67},
\cite[p.~28]{BB}.  
\end{quotation}
The fact that an~$a$ satisfying~$\neg ((a\leq 0) \vee (a \geq 0))$
yields a counterexample~$f(x)=ax$ to the extreme value theorem (EVT)
on~$[0,1]$ is alluded to by Bishop in \cite[p.~59, exercise~9]{Bi67};
\cite[p.~62, exercise~11]{BB}.  Meanwhile, Bridges interprets Bishop's
italicized ``{\em not\/}'' as referring to a Brouwerian
counterexample, and asserts that trichotomy as well as the
principle~$(a\leq 0) \vee (a \geq 0)$ are independent of Bishopian
constructivism.  See D.~Bridges~\cite{Br94} for details; a useful
summary may be found in Taylor~\cite{Ta}.

The extreme value theorem illustrates well the fact that a
verificational interpretation of the quantifiers in Bishopian
mathematics necessarily results in a clash with classical mathematics;
it is merely a matter of tactical emphasis that the clash is minimized
in Bishopian constructivism.

\section{Are there two constructivisms?}

\subsection{Numerical constructivism or anti-LEM constructivism?}
\label{51}

G.~Hellman has distinguished between {\em liberal\/} and {\em
radical\/} versions of constructivism in \cite{Hel93a}.  The
dichotomy%
\footnote{\label{radical}While such a dichotomy seems to have taken
root among the philosophers of mathematics, it may take time to
penetrate the radical defenses.}
was picked up both by by H.~Billinge~\cite{Bil03} and by
E.~Davies~\cite{Da05}.  Rather than borrowing terms from political
science, we will exploit terms that are somewhat more self-explanatory
in a mathematical context.

What are the goals of constructive mathematics?  Bishop in his 1975
{\em Crisis\/} essay \cite{Bi75} emphasizes a quest for {\em numerical
meaning\/} (see Subsection~\ref{eb}).

It seems reasonable to assume that a mainstream mathematician can
relate favorably to the methodological approach that seeks greater
numerical meaning, as enunciated by Bishop in \cite{Bi75}.  Such an
approach could be described as {\em numerical constructivism\/}, or
rather, numerically meaningful constructivism.  A methodological
approach, complementing other approaches, could be a {\em companion\/}
to classical mathematics, while at the same time recognizing the
coherence of post-LEM numerical meaning (see
Subsection~\ref{postlem}).

On the other hand, when the emphasis shifts to the extirpation of the
law of excluded middle from the mathematical toolkit, whether there is
numerical benefit or loss%
\footnote{Specifically, one loses the intuitively appealing
infinitesimal definitions of the calculus, certain results of
mathematical physics, the calculus of variations, soap films and
bubbles (see Subsection~\ref{fifteen}), and other types of post-LEM
numerical meaning discussed in Subsection~\ref{postlem}.}
in such a quest, one arrives at an approach that could be termed {\em
anti-LEM constructivism\/}, conceived of as an {\em alternative\/} to
classical mathematics, and viewing LEM as a variety of mathematical
phlogiston.%
\footnote{See Subsection~\ref{pour} on Pourciau's phlogiston
metaphor.}

To summarize, Bishop's comments could lead to two distinct conceptions
of constructivism:

\begin{enumerate}
\item
numerical constructivism, a {\em companion\/} to classical
mathematics; and
\item
anti-LEM constructivism, an {\em alternative\/} to classical
mathematics.
\end{enumerate}

Returning to Bishop's {\em Crisis\/} essay \cite{Bi75}, we note that
when he found himself in front of a distinguished audience at the
Boston workshop,%
\footnote{\label{distinguish}His audience included Birkhoff,
Dieudonn\'e, Freudental, J.P.Kahane, Kline, and Mackey.}
he chose to emphasize the goal of seeking greater {\em numerical
meaning\/} over the goal of the elimination of LEM.  Such an attitude
is apparently consistent with Heyting's approach discussed above.

Elsewhere, however, Bishop adopts a different approach, implying that
unless one extirpates LEM, one's mathematics falls in the range
between {\em obfuscation\/} and {\em debasement of meaning\/}.  In
short, constructive mathematics is posited, not as a {\em
companion\/}, but rather as a {\em replacement for\/} and {\em
alternative to\/}, classical mathematics.

Should the gentle reader feel that our description of an anti-LEM
species of constructivism is exaggerated, she should ponder
H.~Billinge's comment reproduced below.  Billinge's sympathies for
constructivism have led her to issue a philosophical challenge, so as
\begin{quote}
to settle the matter of whether we should {\bf permit or prohibit} the
use of classical mathematics%
\footnote{\label{doom}Note that in Hellman's view \cite[p.~439]{He98},
``any [\dots] attempt to reinstate a `first philosophical' theory of
meaning prior to all science is doomed''.  What this appears to mean
is that, while there can certainly be a philosophical notion of
meaning before science, any attempt to {\em prescribe\/} standards of
meaning {\em prior\/} to the actual practice of science, is {\em
doomed\/}.  See also discussion of Dummett in main text at
footnote~\ref{f39}.}
\cite[p.~317]{Bil} [emphasis added--authors].
\end{quote}

Elsewhere, Billinge has claimed that E.~Bishop is {\em not\/} a
radical constructivist:
\begin{quote}
We can safely say that Bishop was a liberal constructivist
\cite[p.~187]{Bil03}.
\end{quote}
Yet, one wonders what Billinge would make of Bishop's published
opinion to the effect that
\begin{quote}
Very possibly classical mathematics will cease to exist as an
independent discipline \cite[p.~54]{Bi68}.%
\footnote{\label{all3}See also footnotes~\ref{all1}, \ref{all2a},
and~\ref{all2}.}
\end{quote}
Furthermore, Bishop's vitriolic book review \cite{Bi77} is absent from
Billinge's bibliography, and she makes no attempt to explain Bishop's
``debasement'' comments \cite[p.~513-514]{Bi75} and \cite[p.~1]{Bi85}
targeting classical mathematics as a whole.

\subsection{Hilbert and post-LEM meaning}
\label{postlem}

To elaborate further on the difference between the two approaches to
constructivism discussed in Subsection~\ref{51}, we will examine the
attitude toward calculus, from the perspective of the
formalist/intuitionist debate.  The infinitesimal definition of the
derivative was already envisioned by the founders of the calculus, and
was justified by A.~Robinson (see Appendix~\ref{rival}).  Heyting is
able to appreciate such an accomplishment (see Subsection~\ref{APR}),
firmly grounded, as it was, in classical logic, making it unacceptable
to an anti-LEM constructivist.

An infinitesimal definition certainly leads to computationally
meaningful formulas for the derivatives of the standard functions,
even though classical logic was relied upon by Robinson.  To emphasize
the difference between the two types of numerical meaning, we
introduce the following distinction.  

We will refer to the content of mathematics that is fully constructive
from the bottom (i.e. integers) up, as {\em pre-LEM numerical
meaning\/}, and to the content of mathematics that has indisputable
numerical or computational meaning, while at the same time relying on
LEM at some intermediate stage, as {\em post-LEM numerical
meaning\/}.%
\footnote{\label{rich9}See main text preceding footnote~\ref{richman}
for an analysis of Richman's views in the light of our pre/post-LEM
dichotomy.}

Pre-LEM numerical meaning is, of course, exemplified by {\em
Foundations of constructive analysis\/} \cite{Bi67}.  We will now
present several examples of post-LEM numerical meaning, following
Euclid, Leibniz, Hilbert, Loewner, Connes, Yau, and Avigad.

\begin{example}
An early example of what we consider post-LEM numerical meaning is
given by geometric constructions in antiquity.

The volume calculations in Archimedes appear to rely on proofs by
contradiction.  Thus, a contradiction is derived from the assumption
that the volume~$x$ of the unit sphere is {\em smaller than\/}
two-thirds of the volume of the circumscribed cylinder, and similarly
for {\em greater than\/}.  The contradiction proves that the volume is
exactly {\em two-thirds\/}.  The proof appears to depend on the
trichotomy for real numbers:
\[
(x<y) \vee (x=y) \vee (x>y),
\]
where~$y$ is two-thirds of the volume of the cylinder.  We start with
the trichotomy, eliminate the possibilities on the right and on the
left, to conclude that the possibility in the middle (equality) is the
correct one.  

However, the proof can be modified slightly to meet the requirements
of Bishop's framework.  Thus, proving first that $\neg(x<y)$
implies~$x\geq y$ (but not {\em vice versa\/})%
\footnote{One cannot derive, constructively,~$x<y$ by getting a
contradiction from the assumption that~$x\geq y$.  Note that the
formula~$\neg(x<y)\to (x\geq y)$ implies its contrapositive
$\neg(x\geq y) \to \neg\neg(x<y)$, but double negation~$\neg\neg(x<y)$
does {\em not\/} imply~$x<y$ if one is working constructively.}
and~$\neg(x>y)$ implies~$x\leq y$, to conclude that~$x=y$ without
appealing to trichotomy.

Meanwhile, Euclid's constructions were recently examined, and found to
be lacking, from a constructivist viewpoint, by
M.~Beeson~\cite{Bee09}.  Beeson succeeds in reformulating much of
Euclid in a constructively acceptable fashion, namely without
``test-for-equality'' constructions (such constructions rely on LEM so
as to be able to consider the two cases~$A=B$ and~$A\not=B$
separately).  

Beeson's work is no doubt a significant accomplishment.  Yet one can
question the ideological assumption that, due to their reliance on
LEM, Euclid's arguments were somehow lacking in meaning prior to
Beeson's arrival on the scene.  Will such an assumption be shared by
all of Beeson's readers?
\end{example}

\begin{example}
\label{112}
Leibniz's calculus was based on an apparently computationally
questionable entity.  Indeed, the theoretical entity called the {\em
infinitesimal\/} appeared to have no empirical or computational
counterpart.%
\footnote{\label{sherry}This aspect of the infinitesimal was a source
of Berkeley's metaphysical criticism (Berkeley's criticism was
dissected into its metaphysical and logical components by D.~Sherry
\cite{She87}).  Berkeley's empiricist philosophy stemming from his
perception-based theory of knowledge, tolerated no theoretical
entities without an empirical referent.}
Thus, the algorithms for computing derivatives as developed by
Leibniz, appear to rely upon non-constructive foundational material.
Infinitesimal calculus remains in the post-LEM category even after
Robinson's work, due to the nature of the construction of the
hyperreal number system (but see Palmgren \cite{Pa}).
\end{example}

\begin{example}
\label{action}
The Einstein-Hilbert action (alternatively, Einstein-Hilbert
Lagrangian) is the action that yields Einstein's field equations for
the spacetime metric in general relativity.  S.~P.~Novikov \cite{No2}
writes that, for all of Hilbert's interest in a program of
formalisation and axiomatisation, he understood such a program in a
highly nontrivial fashion.  Novikov mentions, notably, Hilbert's deep
theorem \cite{Hi15} to the effect that the Einstein equations of
relativistic gravitation are of Lagrangian type.  The equations are
thereby placed in the context of a variational problem.  Novikov's
essay is discussed further in Section~\ref{fifteen}, in the context of
%
%
the general difficulty of accounting constructively for variational
principles.%
\footnote{With hindsight, it can be said that a critique of
intuitionism based on variational principles, rather than boxer's
fists, may have been a more effective critique by Hilbert.}
This is only the most famous example of a meaningful variational
problem that possesses an undeniable meaning of a post-LEM variety.
\end{example}

\begin{example}
\label{loewner}
Loewner's torus inequality (see Katz~\cite{SGT}) is proved using the
uniformisation theorem, which is not fully constructive.  Yet the
inequality itself possesses indisputable numerical meaning.  The
inequality relates a pair of metric invariants of a Riemannian
2-torus~$\T$.  The first is the least length, denoted~$\sys$, of a
loop on $\T$ that cannot be contracted to a point on $\T$.  The second
invariant is the total area, denoted~$\area$, of $\T$.  The
inequality, $\sys^2 \leq \tfrac{2}{\sqrt{3}} \area$, is illustrated in
Figure~\ref{micro}.  Mikhail~Gromov writes as follows:
\begin{quote}
I was exposed to metric inequalities in the late 60's by Yu.~Burago,
who acquainted me with the results of Loewner, Pu, and Besicovitch.
These attracted me by the topological purity of their underlying
assumptions, and I was naturally tempted to prove similar inequalities
in a more general topological framework \cite[p.~271-272]{Gr4}.
\end{quote}
It is not even clear that a weaker version~$\sys^2 <
(\tfrac{2}{\sqrt{3}}+\epsilon) \area$ can be obtained in the absence
of LEM.
\end{example}

\begin{example}
\label{alain}
Alain Connes' theory of infinitesimals is based on non-constructive
foundational material.%
\footnote{See footnote~\ref{connes}.}
Connes wrote as follows:
\begin{quote}
Our theory of infinitesimal variables [\dots] will give a precise
computable answer \cite[p.~6207]{Co95}.
\end{quote}
Connes claims that his infinitesimals provide a computable answer to a
probability problem outlined earlier in his text.  A constructivist
who cares to analyze Connes' remarks, is faced with a stark choice:
either accept the coherence of post-LEM numerical meaning as
exemplified by applications of Connes' infinitesimals, or consider the
founder of one of the most active fields of mathematical physics as
being incoherent when talking about ``computable'' answers.
\end{example}

\begin{example}
The existence of Calabi-Yau metrics \cite{YN}, influential both in
differential geometry and in high-energy physics, appears to be
non-constructive, as it relies on the calculus of variations,
including numerous applications of versions of the extreme value
theorem on several levels in infinite-dimensional spaces.  In this
sense, Calabi-Yau metrics pose the same problem to an anti-LEM
constructivist as the Hilbert-Einstein action (see
Example~\ref{action}).
\end{example}

\begin{example}
\label{avigad}
A detailed study by J. Avigad \cite{Av} of ergodic theory can be
viewed as an analysis of post-LEM numerical meaning, in the context of
ergodic theory.  More generally, the technique of {\em proof mining\/}
seeks to extract numerical content from non-constructive proofs, using
logical tools, see U.~Kohlenbach and P.~Oliva~\cite{KO}.
\end{example}

Foundationally speaking, the acceptance of the validity of post-LEM
numerical meaning may well be a fall-back to Hilbert's position,
called {\em syntactic modeling\/} by Avigad \cite{Av}.  It is deemed
acceptable to introduce non-constructive techniques, so long as one
could not derive anything unacceptable from them, where
``unacceptable" means numerically testable results that are found to
be false.%
\footnote{Syntactic modeling may have been at the origin of Hilbert's
interest in proving consistency.  It is generally considered that
Hilbert's hopes for a consistency proof have been dashed by Goedel's
incompleteness theorems, but see Feferman \cite{Fe93} and Avigad and
Reck \cite{AR} for a more nuanced picture.}

Attempting to analyze Hilbert's position, Stolzenberg writes as
follows:

\begin{quote}
Hilbert saw formal mathematics as a way of reaching the real by
passing through the ideal.  That is, one may use all the formal
machinery, in particular, nonconstructive but formally valid existence
statements (such as the Bolzano-Weierstrass theorem), to prove,
formally, real propositions, i.e. predictive ones \cite{Sto}.
\end{quote}
Stolzenberg concludes on the following optimistic note:

\begin{quote}
However, not only are such considerations largely ignored or blurred
nowadays, but we have already quoted, and can confirm, Bishop's
observation that in practice such proofs are already constructive or
can easily be made so \cite{Sto}.
\end{quote}
Stolzenberg examines Hilbert's defense of classical mathematics (based
as it is on classical logic, including LEM), and finds it superfluous.
In our terminology introduced earlier, Stolzenberg appears to claim
that post-LEM numerical meaning ``can easily be made'' pre-LEM.  Now
the constructivisation of classical results may succeed (at any rate
modulo strengthening the hypotheses or weakening the conclusions, see
Subsection~\ref{shesh}) in the case of statements such as the
intermediate value theorem (perhaps even, as implied by Stolzenberg,
the Bolzano-Weierstrass theorem), and the like.  But Bishop apparently
believed that a Robinson infinitesimal resists constructivisation.

\subsection{How could an intuitionist accept the hyperreals?}
\label{53}

If both forms of constructivism discussed in Subsection~\ref{51} are
opposed to the law of excluded middle, then wouldn't Robinson's theory
run afoul of both, seeing that non-standard analysis rests on
non-constructive foundations of classical logic?

To answer the question, we will start by pointing out that the burden
of explanation is upon the constructivists on at least the following
two counts:
\begin{enumerate}
\item
How is it that one of their own, A. Heyting, described Robinson's
theory as not only meaningful but {\em a standard model of important
mathematical research\/}? (see Subsection~\ref{APR} for an analysis of
Heyting's comments).%
\footnote{Heyting gave the first formal development of intuitionistic
logic in order to codify Brouwer's way of doing mathematics; yet, he
unequivocally sides with the view of intuitionism as a {\em
companion\/}, rather than {\em alternative\/}, to classical
mathematics (see Subsection~\ref{APR}).}
\item
to the extent that non-standard analysis shares its foundations with,
and is therefore part of, mainstream mathematics, a rejection of the
former entails a radical rejectionist stance (on the part of the
constructivist) that in particular encounters serious difficulties
{\em vis-a-vis\/} important results of relativity theory, such as the
Hawking--Penrose theorem, see Subsections~\ref{fifteen} and~\ref{HPT}.
\end{enumerate}

Bishop sometimes seems to identify his notion of {\em numerical
meaning\/}, with {\em meaning\/} in a lofty epistemological sense, see
Section~\ref{fourteen}, but Heyting appears to view numerical meaning
as a desirable, but not exclusive, methodological goal.  Thus, Bishop
and Heyting disagree as to the coherence of post-LEM numerical meaning
(see Subsection~\ref{postlem}).

Heyting's position is not inconsistent with what we described in
Subsection~\ref{51} as {\em numerical constructivism\/}, see also
Subsection~\ref{eight} below.  One can adopt the methodological goals
of a search for greater {\em numerical meaning\/}, without locking
oneself into an anti-LEM constructivism.%
\footnote{See Crowe's comment in the main text following
footnote~\ref{crowe}.}

In Heyting's book \cite{He}, a protagonist named {\em Class\/}
(apparently alluding to {\em classical logic\/}) asserts that
``intuitionism should be studied as a part of mathematics''.  Note
that a protagonist named {\em Int\/} (apparently alluding to {\em
intuitionism\/}) readily agrees to {\em Class\/}'s contention (see
Subsection~\ref{eight}).  Is {\em Int\/} merely being polite in
conceding a point to {\em Class\/}?  There is no evidence in Heyting's
book that he is merely being indirect.  On the other hand, one could
well analyze Bishop's comment in this light.  When Bishop writes
\cite{Bi75} that he ``does not know" whether non-standard analysis
results in a loss of meaning, and in the same paragraph proceeds to
characterize non-standard calculus as a {\em debasement\/} of the
latter (see Section~\ref{eb}), one could feel that his comments can be
described as {\em indirect\/}.

\section{The fervor of Bishopian constructivism}
\label{fervor}

In this section, we seek to illustrate Bishop's influence on the
standards of the debate concerning the nature of mathematics.

\subsection{Pourciau on phlogiston}
\label{pour}

The influence of Bishop's radical philosophy has been felt in the area
of discourse and standards of civility, in publications that have
appeared in periodicals of the history and the philosophy of
mathematics.  One of the most rousing affirmations of the anti-LEM
credo in English was penned by B. Pourciau and published in the {\em
Studies in History and Philosophy of Science\/}:
\begin{quote}
The faith that sustains the classical world view emanates from one
belief more than any other: that mathematical assertions are true or
false, independently of our knowing which. Every {\bf conflict}
between classical and intuitionist mathematics springs ultimately from
this belief.  This is not a belief so much as a hidden cause: it
creates the world view---making, for example, proofs by contradiction
appear self-evidently correct---but remains transparently in the
background, unseen and {\bf unquestioned} \cite[p.~319]{Po00}
[emphasis added--authors].
\end{quote}

Such an {\em unquestioning\/} stance naturally leads to {\em
calcification\/}, as Pourciau continues: 
\begin{quote}
Once created, the classical world view is sustained and {\bf
calcified} by metaphor taken literally, by what Stolzenberg%
\footnote{See \cite[p.~225]{Sto78a}.}
calls a 
\begin{quote}
``present tense'' language of ``objects and their properties''.
\end{quote}
Talk about sets, functions, real numbers, theorems and so on is taken
by the classical mathematician as being literally about mathematical
objects that exist independently of us, even though such talk,
classically interpreted, has the appearance, and nothing more, of
being {\bf meaningful}%
\footnote{Pourciau's analysis of classical mathematics conflates two
distinct positions, described by W.~Tait \cite{Ta01} as {\em
realism\/} and {\em super-realism\/}.}
\cite[p.~319]{Po00} [emphasis added--authors].
\end{quote}

The reader recognizes the pet constructivist term, {\em meaning\/},
already analyzed in Subsection~\ref{eb} (see also
Subsection~\ref{glossary} below), as Pourciau continues by quoting
Stolzenberg as follows: `` `To anyone who starts off in the
contemporary mathematician's belief system', argues Stolzenberg
\cite[p. 268]{Sto78a},

\begin{quote}
the discovery that an entire component of the `reality' of one's
experience is produced by {\bf acts of acceptance} as such in the
domain of language use is not merely illuminating. In a literal sense,
it is {\bf shattering}.  Once a mathematician has seen that his
perception of the `self-evident correctness' of the law of excluded
middle is nothing more than the linguistic equivalent of an optical
illusion, neither his practice of mathematics nor his understanding of
it can ever be the same'' \cite[p.~319]{Po00} [emphasis
added--authors].
\end{quote}

Pourciau's introduction, with its talk of ``conflict'', ``acts of
acceptance'', and ``shattering'' discoveries, is a good sample of an
insurrectional narrative, two types of which are analyzed below in
Subsections~\ref{insurrection} and~\ref{insurrection2}.  J.~Avigad,
while agreeing that
\begin{quote}
[w]e do not need fairy tales about numbers and triangles prancing
about in the realm of the abstracta[,]
\end{quote}
notes that 
\begin{quote}
[p]roof-theoretic analysis [\dots] yields satisfying philosophical
explanations as to how abstract, infinitary assumptions have bearing
on computational concerns, and provides senses in which infinitary
methods can be seen to have finitary content \cite{Av07}.%
\footnote{For more details see Example~\ref{avigad} in
Subsection~\ref{postlem} above.}
\end{quote}

Having rounded off his ``shattering'' introduction, Pourciau now
prepares his main thesis by quoting G.~Giorello as follows.

\begin{quote}
One may ask whether there is a `phlogiston'%
\footnote{Phlogiston was once thought to be a fire-like element
contained within combustible bodies, and released during combustion,
but became incommensurable 250 years ago.}
in mathematics. \dots I would be inclined to say `No.'  This, in our
opinion, would constitute a difference between a mathematical
revolution and a `revolution' in Kuhn's sense \cite[p. 168]{Gi}.
\end{quote}

Pourciau takes this as evidence that Giorello does {\em not\/} argue
against the position, attributed to Crowe, that Kuhnian revolutions
are ``inherently impossible'' in mathematics.  Meanwhile, Pourciau
himself argues that, had Brouwer's rejection of LEM been generally
accepted, it would in fact have constituted such a revolution.  Thus
Pourciau appears to understand Intuitionism as viewing LEM as just
such a variety of phlogiston, a position we would characterize as
anti-LEM (see Subsection~\ref{51}).

To elaborate, we argue as follows.  Had Pourciau not thought of LEM as
mathematical phlogiston, he would not have taken Giorello's remarks as
evidence that Giorello does not argue against Crowe's contention that
revolutions never occur in mathematics.  The contradiction proves that
Pourciau does in fact think of LEM as mathematical phlogiston.  Our
argument by contradiction, of course, is only meaningful classically,
and would presumably be incoherent in Pourciau's preferred
intuitionistic framework.  However, we include it here for the benefit
of a classically minded reader.

Pourciau's novel view is that Brouwer's Intuitionism constituted a,
{\em failed,\/} Kuhnian revolution.  Pourciau motivates his view by
dint of the observation that certain classically true assertions
become ``incommensurable'' in Brouwer's framework.  Note that the
essence of Pourciau's view had already been captured by Giorello, who
wrote that

\begin{quote}
``even rigorists like Weierstrass were labelled `metaphysicians' by
strict constructivists like Weyl or Brouwer'' \cite[p.~166]{Gi}.
\end{quote}

Pourciau's main example is the (classically true) assertion to the
effect that

\begin{quote}
in the infinite decimal expansion of~$\pi$, either finitely many pairs
of consecutive equal digits occur, or infinitely many pairs occur
\cite[p.~305]{Po00}.
\end{quote}
Pourciau describes this assertion as becoming ``incommensurable'' in
an intuitionistic framework.  Based on the loss of such assertions in
an intuitionistic framework, Pourciau concludes that Intuitionism
possessed a potential of developing into a unique mathematical Kuhnian
revolution.  Pourciau writes that

\begin{quote}
from inside the classical paradigm [\dots] Kuhnian revolutions {\em
appear\/} to be [\dots] obviously impossible.  But from {\em
outside\/} the classical paradigm, [\dots] it is no longer
self-evident that Kuhnian revolutions cannot occur
\cite[p.~300]{Po00}.
\end{quote}

He argues that both Crowe and Dauben 
\begin{quote}
{\em see Kuhnian revolutions as logically impossible, because to them
[\dots] any shift in conceptions of mathematics must be cumulative.\/}
\cite[p.~301]{Po00} [emphasis in the original--the authors].
\end{quote}
Pourciau views such cumulativity as a necessary consequence of the
classical (i.e. LEM-circumscribed) paradigm, suggesting that a Kuhnian
revolution in mathematics would in fact be impossible, without first
extirpating LEM.

However, contrary to Pourciau's claim, mathematical entities have
indeed become ``incommensurable'' through historical changes in
mathematical practice, ever before intuitionism came on the scene.

Consider, for example, the largely successful assault on
infinitesimals in the aftermath of the rise of Weierstrassian
epsilontics.  Cauchy's 1821 sum theorem \cite{Ca21} on the convergence
(to a continuous limit) of a series of continuous functions under a
suitable ``pointwise'' convergence, was declared to be false in the
absence of infinitesimals (see below).  With hindsight, we can now
affirm that Cauchy's sum theorem had been made ``incommensurable'' in
the second half of the 19th century.  The sum theorem was only
successfully ``revived'' fully a century later.  A similar observation
applies to infinitesimals themselves (see Example~\ref{112}).  Thus,
Weierstrassian epsilontics constitute a, {\em successful\/}, Kuhnian
revolution, by Pourciau's standard.  Indeed, Robinson
\cite[p.~271-273]{Ro66} proposed an interpretation of Cauchy's sum
theorem that would make it correct, in the context of an
infinitesimal-enriched continuum, see also J.~Cleave \cite{Cl},
Cutland {\em et al.\/}~\cite{CKKR}, Goldblatt~\cite[p.~90]{Go}, and
Br\aa ting~\cite{Br07}.  This (ongoing) project appears to be a
striking realisation of a reconstruction project enunciated by
I.~Grattan-Guinness, in the name of Freudenthal \cite{Fr}:
\begin{quote}
it is mere feedback-style ahistory to read Cauchy (and contemporaries
such as Bernard Bolzano) as if they had read Weierstrass already.  On
the contrary, their own pre-Weierstrassian muddles need historical
reconstruction \cite[p.176]{Gra04}.%
\footnote{See Borovik \& Katz \cite{BK} and Katz \& Katz \cite{KK11b}
for a detailed analysis of trends in Cauchy historiography.}
\end{quote}

\subsection{A constructive glossary}
\label{glossary}

As we already mentioned, Bishop's approach is rooted in Brouwer's
revolt against the non-constructive nature of mathematics as practiced
by his contemporaries.%
\footnote{Similar tendencies on the part of Wittgenstein were analyzed
by H.~Putnam, who describes them as ``minimalist''
\cite[p.~242]{Pu07}.  See also G.~Kreisel \cite{Kr}.}

The Brouwer--Hilbert debate captured the popular mathematical
imagination in the 1920s.  Brouwer's crying call was for the
elimination of most of the applications of LEM from meaningful
mathematical discourse.  Burgess discusses the debate briefly in his
treatment of nominalism in \cite[p.~27]{Bu04}.  We analyzed Bishop's
implementation of Brouwer's nominalistic project in~\cite{KK11a}.

Bishop's program has met with a certain amount of success, and
attracted a number of followers.  Part of the attraction stems from a
detailed lexicon developed by Bishop so as to challenge received
(classical) views on the nature of mathematics.  A constructive
lexicon was a {\em sine qua non\/} of his success.  A number of terms
from Bishop's constructivist lexicon constitute a novelty as far as
intuitionism is concerned, and are not necessarily familiar even to
someone knowledgeable about intuitionism {\em per se\/}.  It may be
helpful to provide a summary of such terms for easy reference,
arranged alphabetically, as follows.

\medskip
$\bullet$ {\em Debasement of meaning\/} is the cardinal sin of the
classical opposition, from Cantor to Keisler,%
\footnote{\label{keisler1}But see footnote~\ref{keisler2}.}
committed with {\em LEM\/} (see below).  The term occurs in Bishop's
{\em Schizophrenia\/} \cite{Bi85} and {\em Crisis\/}~\cite{Bi75}
texts.

\medskip
$\bullet$ {\em Fundamentalist excluded thirdist\/} is a term that
refers to a classically-trained mathematician who has not yet become
sensitized to implicit use of the law of excluded middle (i.e.,~{\em
excluded third\/}) in his arguments, see \cite[p.~249]{Ri}.%
\footnote{\label{rich1}This use of the term ``fundamentalist excluded
thirdist'' is in a text by Richman, not Bishop.  I have not been able
to source its occurrence in Bishop's writing.  In a similar vein, an
ultrafinitist recently described this writer as a ``choirboy of
infinitesimology''; however, this term does not seem to be in general
use.  See also Subsections~\ref{richman} and \ref{81}.}

\medskip
$\bullet$ {\em Idealistic mathematics\/} is the output of Platonist
mathematical sensibilities, abetted by a metaphysical faith in {\em
LEM\/} (see below), and characterized by the presence of merely a {\em
peculiar pragmatic content\/} (see below).

\medskip
$\bullet$ {\em Integer\/} is the revealed source of all {\em
meaning\/} (see below), posited as an alternative foundation
displacing both formal logic, axiomatic set theory, and recursive
function theory.  The integers wondrously escape%
\footnote{By dint of a familiar oracular quotation from Kronecker; see
also main text around footnote~\ref{potential}.}
the vigilant scrutiny of a constructivist intelligence determined to
uproot and nip in the bud each and every Platonist fancy of a concept
{\em external\/} to the mathematical mind.

In Bishop's system, the integers are uppermost to the exclusion of the
continuum.  Bishop rejected Brouwer's work on an intuitionstic
continuum in the following terms:
\begin{quote}
Brouwer's bugaboo has been compulsive speculation about the nature of
the continuum.  His fear seems to have been that, unless he personally
intervened to prevent it, the continuum would turn out to be discrete.
[The result was Brouwer's] semimystical theory of the continuum
\cite[p.~6 and 10]{Bi67}.
\end{quote}

\medskip
$\bullet$ {\em Integrity\/} is perhaps one of the most misunderstood
terms in Errett Bishop's lexicon.  Pourciau in his {\em
Education\/}~\cite{Po99} appears to interpret it as an indictment of
the ethics of the classical opposition.  Yet in his {\em
Schizophrenia\/} text, Bishop merely muses:
\begin{quote}
I keep coming back to the term ``integrity''.  \cite[p.~4]{Bi85}
\end{quote}
Note that the period is in the original.  Bishop describes {\em
integrity\/} as the opposite of a syndrome he colorfully refers to as
{\em schizophrenia\/}, characterized by a number of ills, including
\begin{enumerate}
\item[(a)] a rejection of common sense in favor of formalism, 
\item[(b)] the {\em debasement of meaning\/} (see above),
\item[(c)] 
as well as by a list of other ills---
\end{enumerate}
but {\em excluding\/} dishonesty.  Now the root of
\begin{equation*}
\mbox{integr-ity}
\end{equation*}
is identical with that of {\em integer\/} (see above), the Bishopian
ultimate foundation of analysis.  Bishop's evocation of {\em
integrity\/} may have been an innocent pun intended to allude to a
healthy constructivist mindset, where the {\em integers\/} are
uppermost.

Brouwer sought to incorporate a theory of the continuum as part of
intuitionistic mathematics, by means of his free choice sequences.
Bishop's commitment to integr-ity is thus a departure from Brouwerian
intuitionism.

\medskip
$\bullet$ {\em Law of excluded middle (LEM)\/} is the main source of
the non-constructivities of classical mathematics.%
\footnote{See Example~\ref{roottwo}.}
Every formalisation of {\em intuitionistic logic\/} excludes {\em
LEM\/}; adding {\em LEM\/} back again returns us to {\em classical
logic\/}.

\medskip
$\bullet$ {\em Limited principle of omniscience (LPO)\/} is a weak
form of {\em LEM\/} (see above), involving {\em LEM\/}-like oracular
abilities limited to the context of integer sequences.%
\footnote{See Subsection~\ref{eb} for a discussion of LPO.}
The {\em LPO\/} is still unacceptable to a constructivist, but could
have served as a basis for a {\em meaningful\/} dialog between Brouwer
and Hilbert (see \cite{Bi75}), that could allegedly have changed the
course of 20th century mathematics.

\medskip
$\bullet$ {\em Meaning\/} is a favorite philosophic term in Bishop's
lexicon, necessarily preceding an investigation of {\em truth\/} in
any coherent discussion.  In Bishop's writing, the term {\em
meaning\/} is routinely conflated with {\em numerical meaning\/} (see
below).

\medskip
$\bullet$ {\em Numerical meaning\/} is the content of a theorem
admitting a proof based on intuitionistic logic, and expressing
computationally meaningful facts about the integers.%
\footnote{As an illustration, a numerically meaningful proof of the
irrationality of~$\sqrt{2}$ appears in Example~\ref{roottwo}.}
The conflation of {\em numerical meaning\/} with {\em meaning\/} par
excellence in Bishop's writing, has the following two consequences:
\begin{enumerate}
\item
it empowers the constructivist to sweep under the rug the distinction
between pre-LEM and post-LEM numerical meaning,%
\footnote{See Subsection~\ref{postlem}.}
lending a marginal degree of plausibility to a dismissal of classical
theorems which otherwise appear eminently coherent and meaningful;%
\footnote{See Subsection~\ref{evt} for a discussion of the extreme
value theorem in an intuitionstic context.}
and
\item
it allows the constructivist to enlist the support of {\em
anti-realist\/} philosophical schools of thought (e.g.~Michael
Dummett) in the theory of meaning, inspite of the apparent tension
with Bishop's otherwise {\em realist\/} declarations (see entry {\em
Realistic mathematics\/} below).
\end{enumerate}

\medskip
$\bullet$ {\em Peculiar pragmatic content\/} is an expression of
Bishop's \cite[p.~viii]{Bi67} that was analyzed by Billinge
\cite[p.~179]{Bil03}.  It connotes an alleged lack of empirical
validity of classical mathematics, when classical results are merely
{\em inference tickets\/} \cite[p.~180]{Bil03} used in the deduction
of other mathematical results.

\medskip
$\bullet$ {\em Realistic mathematics\/}.  The dichotomy of ``realist''
{\em versus\/} ``idealist'' (see above) is the dichotomy of
``constructive'' {\em versus\/} ``classical'' mathematics, in Bishop's
lexicon.  There are two main narratives of the Intuitionist
insurrection, one {\em anti-realist\/} and one {\em realist\/}.  The
issue is discussed in Subsections~\ref{insurrection} and
\ref{insurrection2}.

On the first page of his book, Bishop claims that constructive
mathematics based on the natural numbers will still be the {\em
same\/}, ``in another universe, with another biology and another
physics'' \cite[p.~1]{Bi67}.  Since the sensory perceptions of the
human {\em body\/} are physics- and chemistry-bound, a claim of such
trans-universe invariance amounts to the positing of a {\em
disembodied\/} nature of the natural number system (transcending the
physics and the chemistry).  Bishop's disembodied natural numbers are
the cornerstone of his approach.
%

\subsection{Insurrection according to Michael}
\label{insurrection}

The {\em anti-realist\/} narrative, mainly following Michael Dummett
\cite{Du}, traces the original sin of classical mathematics with {\em
LEM\/}, all the way back to Aristotle.%
\footnote{\label{keisler2}the entry under {\em debasement of
meaning\/} in Subsection~\ref{glossary} would read, accordingly, ``the
classical opposition from Aristotle to Keisler''; see main text at
footnote~\ref{keisler1}.}
The law of excluded middle (see Subsection~\ref{glossary}) is the
mathematical counterpart of geocentric cosmology (alternatively, of
phlogiston, see Subsection~\ref{pour}), slated for the dustbin of
history.%
\footnote{\label{dummett}Following Kronecker and Brouwer, Dummett
rejects actual infinity, at variance with Bishop.}
The anti-realist narrative dismisses the Quine-Putnam indispesability
thesis (see Subsection~\ref{fifteen}), on the grounds that a {\em
philosophy-first\/} examination of first principles is the unique
authority empowered to determine the correct way of doing
mathematics.%
\footnote{\label{f39}See discussion of the views of Billinge and
Hellman in the main text at footnote~\ref{doom}.}
Generally speaking, it is this narrative that seems to be favored by a
number of philosophers of mathematics.

Dummett opposes a truth-valued, bivalent semantics, namely the notion
that truth is one thing and knowability another, on the grounds that
it violates Dummett's {\em manifestation requirement\/}, see
Shapiro~\cite[p.~54]{Sh}.  The latter requirement, in the context of
mathematics, is a restatement of the intuitionistic principle that
truth is tantamount to verifiability.  Thus, an acceptance of
Dummett's manifestation requirement, leads to intuitionistic semantics
and a rejection of LEM.

In his foundational text \cite{Du77} originating from 1973 lecture
notes, Dummett is frank about the source of his interest in the
intuitionist/classical dispute in mathematics:
\begin{quote} 
This dispute bears a {\bf strong resemblance} to other disputes over
realism of one kind or another, that is, concerning various kinds of
subject-matter (or types of statement), including that over realism
about the physical universe%
\footnote{The {\em strong resemblance\/} claim is precisely the point
we wish to dispute, see below.}
\cite[p.~ix]{Du77} [emphasis added--authors].
\end{quote}
What Dummett proceeds to say at this point, reveals the nature of his
interest:
\begin{quote} 
but intuitionism represents the only sustained attempt by the
opponents of a {\bf realist view} to work out a coherent embodiment of
their philosophical beliefs [emphasis added--authors].
\end{quote}
What interests Dummett here is the fight against the {\em realist
view\/}.  What endears intuitionists to him, is the fact that they
have succeeded where the phenomenalists have not:
\begin{quote}
Phenomenalists might have attained a greater success if they had made
a remotely comparable effort to show in detail what consequences their
interpretation of material-object statements would have for our
employment of our language \cite[p.~ix]{Du77}.
\end{quote}
However, Dummett's conflation of the mathematical debate and the
philosophical debate, could be challenged.

We hereby explicitly sidestep the debate opposing the realist (as
opposed to the super-realist, see W.~Tait \cite{Ta01}) position and
the anti-realist position.  On the other hand, we observe that a
defense of indispensability of mathematics would necessarily start by
challenging Dummett's manifestation requirement.  More precisely, such
a defense would have to start by challenging the extension of
Dummett's manifestation requirement, from the realm of philosophy to
the realm of mathematics.  While Dummett chooses to pin the opposition
to intuitionism, to a belief in an
\begin{quote}
interpretation of mathematical statements as referring to an
independently existing and objective reality \cite[p.~ix]{Du77},
\end{quote}
(i.e. a Platonic world of mathematical entities), J.~Avigad memorably
retorts as follows:
\begin{quote}
We do not need fairy tales about numbers and triangles prancing about
in the realm of the abstracta ~\cite{Av07}.
\end{quote}
A.~Weir writes that Dummett's argument hinges on questionable
empiricist dogma:
\begin{quote} 
Dummett claims that no teacher could manifest and communicate to a
pupil, about to acquire a language, grasp of sentences which are
determinately true or false though we may be incapable of determining
which.  But the argument seems to hinge on the empiricist dogma that
what is learned cannot radically transcend the experiential or
stimulatory inputs to learning \cite[p.~465]{We86}.
\end{quote}

\subsection{Insurrection according to Errett}
\label{insurrection2}

Turning now to the {\em realist\/} narrative of the intuitionist
insurrection, we note that such a narrative appears to be more
consistent with what Bishop himself actually wrote.  In his
foundational essay, Bishop expresses his position as follows:
\begin{quote}
As pure mathematicians, we must decide whether we are playing a game,
or whether our theorems describe an external reality.%
\footnote{This remark of Bishop's has been reflected in secondary
sources; see \cite[p.~188]{Da88}.}
\cite[p.~507]{Bi75}.
\end{quote}
The right answer, to Bishop, is that they do describe an external
reality.%
\footnote{Our purpose here is not to endorse or refute Bishop's views
on this point, but rather to document his actual position, which
appears to diverge from Dummett's.}
The dichotomy of ``realist'' {\em versus\/} ``idealist'' is the
dichotomy of ``constructive'' {\em versus\/} ``classical''
mathematics, in Bishop's lexicon (see entry under {\em idealistic
mathematics\/} in Section~\ref{glossary}).  In~\cite[p.~4]{Bi85},
Bishop mentions his ambition to incorporate {\em such mathematically
oriented disciplines as physics\/} as part of his constructive
revolution, revealing a recognition, on his part, of the potency of
the Quine-Putnam indispensability challenge.%
\footnote{Billinge \cite[p.~314]{Bil} purports to detect ``inchoate''
anti-realist views in Bishop's writings, but provides no constructive
proof of their existence, other than a pair quotes on {\em numerical
meaning\/}.  Meanwhile, Hellman writes: ``Some of Bishop's remarks
(1967) suggest that his position belongs in [the radical] category''
\cite[p.~222]{Hel93a}.}

N.~Kopell and G.~Stolzenberg, close associates of Bishop, published a
three-page {\em Commentary\/} \cite{KS} following Bishop's {\em
Crisis\/} text.  Their note places the original sin with {\em LEM\/}
at around 1870 (rather than Greek antiquity), when the ``flourishing
empirico-inductive tradition'' began to be replaced by the ``strictly
logico-deductive conception of pure mathematics''.%
\footnote{The latter term may be referring to the Fregean revolution
in logic, see footnote~\ref{frege}.}
Kopell and Stolzenberg don't hesitate to compare the {\em
empirico-deductive tradition\/} in mathematics, to physics, in the
following terms:
\begin{quote}
[Mathematical] theories were theories about the phenomena, just as in
a physical theory \cite[p.~519]{KS}.
\end{quote}
Similar views have been expressed by D.~Bridges, see e.g.~\cite{Br99},
and, as we argue in Subsection~\ref{eight}, by Heyting.

\subsection{Richman on large cardinals}
\label{nine}

The {\em Interview with a constructive mathematician\/} was published
by leading constructivist F.~Richman in 1996.  In the {\em
Interview\/}, Richman seems to reject any alternative to an anti-LEM
constructivism, in his very first comment:
\begin{quote}
the constructive mathematician dismisses classical mathematics as an
exercise in formal logic, much like investigating the consequences of
large cardinal axioms \cite[p.~248]{Ri}.
\end{quote}

It should be noted that fellow constructivist D.~Bridges \cite{Br97}
appears to distance himself somewhat from Richman's {\em Interview\/},
which he described in 1997 as
\begin{quote}
a fascinating, provocative, and by no means standard view on
constructivism in mathematics.%
\footnote{Philosophers of science such as Hellman (see
\cite[p.~222]{Hel93a}) and Billinge see Brigdes as being closer to the
liberal outlook, and further from the radical one, than Richman,
whether or not Bridges meant this comment as a way of distancing
himself from Richman.}
\end{quote}
In the context of the exchange between {\em Class\/} and {\em Int\/}
(see Subsection~\ref{eight}), note that Richman seems to concede a
point to {\em Class\/} in the following {\em group theory\/} remark,
somewhat parallel to our {\em Pappus\/} analogy%
\footnote{See discussion at footnote~\ref{pappus}.}
in
Section~\ref{eight}:
\begin{quote}
The richness of constructive mathematics lies in the fact that
concepts that are equivalent in the presence of the law of excluded
middle, need not be equivalent.  This is typical of generalizations:
the notion of a normal subgroup is equivalent to that of a subgroup in
the context of the commutative law \cite[p.~257]{Ri}.
\end{quote}
Yet he quickly retreats to the safety of anti-LEM, observing that
\begin{quote}
the law of excluded middle obliterates the notion of positive content
[\dots]
\end{quote}
Anti-LEM constructivist examination of the foundations is typically
laced with anti-theological sarcasm, as when a mainstream
mathematician is described as ``a fundamentalist excluded-thirdist''
\cite[p.~249]{Ri} (see \cite[p.~260]{Ri} for additional deistic
sarcasm), betraying a nonnegotiable commitment to an ideology.%
\footnote{Whether or not Richman meant this comment to be
tongue-in-cheek, constructivists are merely following Bishop's
sarcastic lead, see \cite[p.~2]{Bi67} as well as Bishop's
``diabolical'' limerick entitled {\em Formalism\/}, characterizing the
latter as ``sawdust'', while constructivism as ``the heart''
\cite[p.~14]{Bi85}.}
Yet we detect a glimmer of hope here, when Richman writes:
\begin{quote}
It now seems apparent to me--although I did not realize this for {\bf
many years}--that for all practical purposes, constructive mathematics
coincides with mathematics done in the context of intuitionistic logic
\cite[p.~253]{Ri} [emphasis added--authors].
\end{quote}
While little evidence is offered for such a counter-intuitive (no pun
intended) assertion, it sets out the following hope.  If a numerical
constructivist can, after {\em many years\/}, be led to abandon
numerical constructivism and switch to anti-LEM; then also, at some
future time, anti-LEM constructivists can perhaps be persuaded, by
force of overwhelming evidence, to abandon anti-LEM and switch to
numerical constructivism of the {\em companion\/} variety.  The nature
of such evidence will be discussed below.

Richman's disillusionment with classical mathematics was triggered by
a perceived lack of meaning in uncountable, torsion abelian group
theory.  Richman writes that, even in the countable case,
\begin{quote}
the centerpiece of the subject, the classification theorem for
countable torsion abelian groups, cannot even be stated without
ordinal numbers \cite{Ri94}.
\end{quote}
To Richman, the non-constructivity of the arguments appears to be
compounded by the non-constructive formulation of the very {\em
statements\/} of the results of torsion abelian group theory.  Richman
apparently perceived a lack of numerical meaning of even the post-LEM
kind, as discussed in Subsection~\ref{postlem}.  If Richman's
disillusionment with classical mathematics was predicated on a
conflation of the two types of numerical meaning, then it leaves open
the possibility for other constructivists to sympathize with Richman's
disillusionment with torsion abelian groups, without, however,
jettisoning post-LEM numerical meaning of classical mathematics.%
\footnote{\label{richman}See also footnote~\ref{rich9}.}

What is the nature of the evidence in favor of a numerical
constructivism of a {\em companion\/} variety?  A mathematician
working in the tradition of Archimedes, Leibniz, and Cauchy owes at
least a residual allegiance to the idea that the most important
mathematical problems are those coming from physics, engineering, and
science more generally.  A mathematical theory that has wide-ranging
applications, be it Abraham Robinson's infinitesimals or Alain Connes'
infinitesimals,%
\footnote{\label{connes}A.~Connes criticizes Robinson's infinitesimals
for being dependent on non-constructive foundational material.  He
further claims it to be a {\em weakness\/} of Robinson's
infinitesimals that the results of calculations that employ them, do
not depend on the choice of the infinitesimal.  Yet, Connes himself
develops a theory of infinitesimals bearing a striking similarity to
the ultrapower construction of the hyperreals.  Furthermore, he freely
relies on such results as the existence of the Dixmier trace, and the
Hahn-Banach theorem.  The latter results rely on similarly
nonconstructive foundational material.  Connes claims the independence
of the choice of Dixmier trace to be a {\em strength\/} of his theory
of infinitesimals in \cite[p.~6213]{Co95}.  Thus, both of Connes'
criticisms apply to his own theory of infinitesimals.  As far as the
introduction of non-commutativity by Connes is concerned, it is
similar to the enlargement from the diagonal matrices, to the spectrum
of arbitrary symmetric matrices.  Such an enlargement is a natural
generalisation.  Its commutative antecedent cannot be criticized for
commutativity, any more than one can criticize the Cartan subalgebra
of, say, a matrix algebra, for being abelian.  The mathematical
novelty of Connes' theory of infinitesimals resides in the
exploitation of Dixmier's trace, relying as it does on
non-constructive foundatoinal material, thus of similar foundational
status to, for instance, the ultrapower construction of a
non-Archimedean extension of the reals.  See also Example~\ref{alain}
above.}
is not easily discounted as meaningless, even if its
starting point may rely on foundational material that is not fully
constructive.  This thread is pursued in more detail in
Section~\ref{fourteen}.

\section{Is constructive mathematics part of classical mathematics?}
\label{partof}

\subsection{From Cantor and Frege to Brouwer}
\label{71}

W.~Tait \cite{Ta83} argues that, unlike intuitionism, constructive
mathematics is part of classical mathematics.  We would like to go
further, and suggest that Frege's revolutionary logic (see Gillies
\cite{Gi92}) and Cantor's revolutionary set theory (see Dauben
\cite{Da92a}) created a new language and a new paradigm, transforming
mathematical foundations into fair game for further investigation,
experimentation, and criticism, including intuitionistic criticism
voiced by Brouwer.%
\footnote{\label{frege}Cantor's set theory is familiar to a general
mathematical audience, and requires no special comment.  The
revolutionary nature of Frege's work from 1879, cf.~\cite{Fre}, is
widely known in logic and philosophy circles.  To summarize Gillies'
thesis, the Fregean revolution in logic was a change from the
Aristotelian paradigm: the theory of syllogism was replaced by
propositional calculus and first-order predicate calculus.  The basis
for the logical notation we today take for granted was in fact laid by
Frege, under the name {\em Begriffsschrift\/}, `conceptual notation'.
(Frege himself used a~$2$-dimensional notation that has not been
accepted by logicians, who have given preference to a linear notation
more along the lines of Peano, who was strongly influenced by Frege.)
For a fascinating account of the genesis of the {\em
Begriffsschrift\/} in Frege's `Logicist program' of showing that the
truths of arithmetic are analytic (rather than synthetic), see
Gillies~\cite[p.~288]{Gi92}.}

Classical mathematics as pursued by Frege and Cantor literally created
a universe of discourse where Brouwer's insurrectional narrative could
be articulated and enjoy a hearing.  In this sense, intuitionism is
not a rival, but an offspring, of classical mathematics.  To quote
J.~Avigad,
\begin{quote}
[t]he syntactic, axiomatic standpoint has enabled us to fashion formal
representations of various foundational stances, and we now have
informative descriptions of the types of reasoning that are justified
on finitist, predicative, constructive, intuitionistic, structuralist,
and classical grounds \cite{Av07}.
\end{quote}
This idea, as applied to intuitionism, was expressed by J.~Gray in the
following terms:
\begin{quote}
Intuitionistic logics were developed incorporating the logical
strictures of Brouwer; constructivist mathematics still enjoys a
certain vogue.  But these are somehow {\bf contained within the larger
framework} of modern mathematics.%
\footnote{Gray concludes with the following thought on current debates
in the foundations of mathematics: ``either they are technical and
accessible only to logicians, or they are epistemological and draw
their examples from concepts we meet in school''.  Examples going
beyond the concepts we meet in school are discussed in
Subsections~\ref{postlem} and~\ref{HPT}.}
\cite[p.~242]{Gray} [emphasis added--authors].
\end{quote}

In a similar vein, no less an authority than M.~Heidegger wrote
(almost simultaneously with Kolmogorov quoted in Subsection~\ref{kol})
as follows:
\begin{quote}
The {\bf level which a science has reached} is determined by how far
it is {\em capable\/} of a crisis in its basic
concepts~\cite[p.~29-30]{Hei} [emphasis added--authors].
\end{quote}
Lest one should doubt whether Heidegger meant for his comments to
apply to mathematics, he continues:

\begin{quote}
{\em Mathematics\/}, which is seemingly the most rigorous and most
firmly constructed of the sciences, has reached a crisis in its
`foundations'.
\end{quote}
Lest one should doubt whether he had Brouwer in mind, Heidegger
continues:

\begin{quote}
In the controversy between the formalists and the {\bf intuitionists},
the issue is one of obtaining and securing the primary way of access
to what are supposedly the objects of this science [emphasis
added--authors].
\end{quote}

What Heidegger appears to be saying is that if we take the supposed
objects to be sets or integers, the issue becomes whether a primary
way of accessing them is provided by the law of excluded middle.  Our
interest in Heidegger's remarks stems from his observation that a
science capable of such a foundational crisis, had necessarily {\em
reached\/} an appropriate {\em level\/} of maturity, which we
attribute to earlier developments due to Frege (via Peano), Cantor,
and others.

\subsection{Heyting's good right}
\label{eight}

A.~Heyting gave the first formal development of intuitionistic logic
in order to codify Brouwer's way of doing mathematics.  Yet, we will
show that Heyting unequivocally sides with the view of intuitionism as
a {\em companion\/}, rather than {\em alternative\/}, to classical
mathematics (see also Subsection~\ref{APR}).  Heyting's example shows
that a concern over the non-constructivity of the widespead use of the
law of excluded middle in classical mathematics, need not necessarily
lead to a radical position expressed in some of Bishop's writing.

A. Heyting's book {\em Intuitionism\/} \cite{He} dates from 1956.  The
(first part of the) book is written in the form of a dialog among
representatives of some of the main schools of thought in the
foundations of mathematics.  We will now analyze Heyting's comments in
more detail.  The protagonist named {\em Class\/} makes the following
point:
\begin{quote}
Intuitionism should be studied as part of mathematics.  In
mathematics, we study consequences of given hypotheses.  The
hypotheses assumed by intuitionists may in fact be interesting, but
they have no right to a {\bf monopoly} \cite[p~.4]{He} [emphasis
added--authors].
\end{quote}
In other words, emptying our logical toolkit of the law of excluded
middle is one possible foundational framework among others.%
\footnote{See footnote~\ref{lem} on the non-uniqueness of
intuitionistic logic.}
This could be compared, for example, to studying projective geometry
with Pappus' axiom eliminated.%
\footnote{\label{pappus}A related point in the context of group theory
is made by F.~Richman \cite{Ri}, see Subsection~\ref{nine}.}
{\em Class\/}'s claim appears to be accepted by the protagonist {\em
Int\/}.  Consider {\em Int's\/} response to {\em Class\/}'s remarks:
\begin{quote}
Nor do we claim that; we [i.e., intuitionists] are content if you
admit the {\bf good right} of our conception [emphasis
added--authors].
\end{quote}
In this connection, S.~Shapiro notes:
\begin{quote}
The later Heyting did not claim a ``monopoly'' on mathematics, and
would rest content if the classical mathematician ``admits the good
right of'' the intuitionistic conception \cite[p.~55-56]{Sh}.
\end{quote}

One can legitimately pose the question whether {\em Int\/} is merely
being polite in his response to {\em Class\/}.  Note that later on in
the same chapter, {\em Int\/} concedes that intuitionism has few
applications, and compares it to ``history, art, and light
entertainment'' \cite[p.~10]{He}.  He does this, obviously, not as a
way of disparaging intuitionism, but as a way of underscoring its
intrinsic value as an intellectual pursuit, as distinct from a
scientific pursuit.  P.~Maddy claims that
\begin{quote}
Heyting [\dots] is happy to dismiss%
\footnote{We saw above that Shapiro clearly disagrees with Maddy on
this point.}
the loss of non-constructive mathematics as `the excision of noxious
ornaments' \cite[pp.~1121-1122]{Ma},
\end{quote}
but a closer examination of Heyting's remark \cite[p.~11]{He} reveals
a more nuanced picture.  Here {\em Int\/} is responding to {\em
Form\/}'s charge to the effect that constructivists
\begin{quote}
destroy the most precious mathematical results \cite[p.~10]{He}.
\end{quote}
{\em Int\/} replies as follows:
\begin{quote}
As to the {\bf mutilation} of mathematics of which you accuse me, it
must be taken as an inevitable consequence of our standpoint.  It can
also be seen as an excision of noxious ornaments, [\dots] and it is at
least {\bf partly compensated for\/} by the charm of subtle
distinctions and witty methods by which intuitionists have enriched
mathematical thought [emphasis added--authors].
\end{quote}
His comments show a clear recognition of the potency of the {\em
mutilation\/} challenge, as well as a necessity to {\em compensate\/}
for the damages. 

Heyting's position is not inconsistent with a {\em numerical
constructivism\/}, as opposed to an anti-LEM constructivism, a kind of
a long march through the foundations launched in {\em sixty-seven\/}
\cite{Bi67}.  Heyting was able to, so to speak, rise above differences
of {\em Class\/} and {\em Int\/}.

What would Heyting have thought of a rigorous justification of
infinitesimals in the framework of classical logic?  Remarkably, we
have a detailed testimony by Heyting himself as to the respect he bore
for Robinson's accomplishment, see Subsection~\ref{APR}.  Already in
1956, Heyting allowed his protagonist {\em Form\/} to cite a reference
by A.~Robinson \cite{Ro51} in a comment dealing with
\begin{quote}
the use of metamathematics for the deduction of mathematical results
\cite[p.~45]{He}.
\end{quote}
In fact, in 1973, Heyting compared Robinson's work favorably with that
of Brouwer himself, according to a Math Reviews summary:
\begin{quote}
[Heyting] gives a brief description of the nature of Robinson's work
in model theory and reflects on how that work resembles that of
Brouwer and thereby qualifies for the award of the Brouwer medal
\cite{He73}.
\end{quote}

S.~Shapiro \cite{Sh}, in a retort to N.~Tennant \cite{Te}, unveils an
allegedly irreducible clash between classical mathematics, on the one
hand, and intuitionism as expressed through Heyting's semantics, on
the other.  Here {\em Heyting semantics\/} refers to the constructive
interpretation of the quantifiers discussed in Subsection~\ref{two}.
Shapiro writes as follows:
\begin{quote}
Let~$\chi(x)$ be any predicate that applies to natural numbers.  It is
a routine theorem of classical arithmetic that
\[
\forall x \exists ! y \left[ (\chi(x) \wedge y=0) \vee (\neg \chi(x)
\wedge y=1) \right] 
\]
(Shapiro \cite[p.~60]{Sh}).
\end{quote}
Classically, Shapiro's formula merely asserts that~$\chi(x)$ is either
true or false (where~$y$ is assigned respectively the value~$0$
or~$1$).  Shapiro continues:
\begin{quote}
Under Heyting semantics, this proof amounts to a thesis that there is
a computable function that decides whether~$\chi$ holds.  So under
Heyting semantics [\dots] every predicate is effectively decidable[,
which is] a tough pull to swallow and keep down.%
\footnote{A similar argument involving the introduction of an
auxiliary bivalent variable occurs in the proof of the Goodman-Myhill
theorem, see footnote~\ref{bm}.}
\end{quote}
Namely, the constructive existence of~$y$ necessarily entails being
able to find a computational procedure capable of identifying such a
$y$ as a function of~$x$, and hence a decision procedure for the
predicate~$\chi$.

Does this prove that classical mathematics and intuitionism {\em
cannot get along\/}, as Shapiro puts it?  {\em Yes\/}, if by
``classical mathematics'' is meant ``all of classical mathematics'',
down to the very last predicate.  On the other hand, the answer is
{\em no\/}, if we are to adopt Bishop's 1975 program~\cite{Bi75},
involving a comprehensive sifting of classical material.

In fact, the philosophical exchange between Tennant and Shapiro has an
uncanny parallel---a quarter century earlier, in the exchange between
Bishop and Timothy G. McCarthy, writing for Math Reviews~\cite{Mc} in
1975.  In the context of Bishop's proposal to translate absolute
classical statements into relative assertions of the form
``LPO$\rightarrow$such-and-such'', McCarthy writes as follows:
\begin{quote}
There is an ambiguity about just what the author intends to explicate
by means of the suggested translation.  It is not plausible to suppose
that~$A$ is classically true (false) if and only if [the proposition]
$\text{LPO}\rightarrow A$ ([respectively,]~$\text{LPO}\rightarrow\neg
A$) is constructively provable, in view of G\"odel's first
incompleteness theorem.
\end{quote}
It appears that Bishop's intention was not to capture all classical
statements, but rather to undertake a comprehensive sift.  Such seems
to be the position adopted by the mature Heyting, as well, as we
analyze in the next section.

\subsection{Heyting's Address to Professor Robinson}
\label{APR}

Intuitionist Heyting's comments in {\em Address to Professor
A. Robinson\/} \cite{He73} are fascinating and worth reproducing in
some detail:
\begin{enumerate}
\item
(p.~135) Brouwer was right that intuitionistic mathematics is the only
form of mathematics which has a perfectly clear interpretation [\dots]
and it is desirable that as much of mathematics as possible will be
made constructive.
\item
(ibid.) the curious mixture of formal reasoning and more or less vague
intuitions which classical mathematics is, remains an imposing work of
art and [at] the same time the powerful tool in the struggle for life
of mankind.%
\footnote{\label{curious}The curious mixture was the subject of
Crowe's musings cited in the main text preceding
footnote~\ref{crowe}.}
\item
(ibid.) intuitionistic mathematics is {\bf no longer isolated} from
classical mathematics [emphasis added--authors].
\item
(ibid.) several mathematicians work on the foundations of
intuitionistic as well as of classical mathematics.
\item
(ibid.) The two subjects become more and more intertwined.
\end{enumerate}
In a remarkable tribute to Robinson's {\em oeuvre\/}, Heyting writes:
\begin{quote}
In your work on the metamathematics of algebra you explored
systematically the connections between formal logic and algebra,
firstly in the form of so-called {\bf transfer principles}, by which a
theorem, proved for some algebraic system, holds automatically for a
whole class of such systems \cite[p.~135]{He73} [emphasis
added--authors].
\end{quote}
For a discussion of the transfer principle of non-standard analysis,
see Appendix~\ref{rival}.  Heyting does not stop there, and makes the
following additional points:
\begin{enumerate}
\item
(p.~136) the creation of non-standard analysis is a standard model of
important mathematical research.
\item
(ibid.) you connected this extremely abstract part of model theory
with a theory apparently so far apart as the {\bf elementary calculus}
[emphasis added--authors].
\item
(ibid.) In doing so you threw new light on the history of the
calculus by giving a clear sense to Leibniz's notion of
infinitesimals.
\item
(ibid.) infinitesimals [\dots] are still commonly used in many fields
where their elimination [\dots] would lead to unwelcome complications.
\item
(ibid.) This is the case in differential geometry and in many parts
of applied mathematics such as hydrodynamics and electricity theory,
where physicists use infinitesimals without a twinge of conscience.
\item
(ibid.) You showed how these theories can be made precise by your
method.  Dirac's~$\delta$-function is an extreme example.  This
mysterious object {\bf became lucid} in the light of non-archimedean
analysis%
\footnote{Infinitesimal constructions of ``Dirac'' delta functions
originate with Fourier and Cauchy, see Borovik \& Katz \cite{BK}.}
[emphasis added--authors].
\end{enumerate}

On the subject of whether there is meaning after model theory, Heyting
has this to say:
\begin{enumerate}
\item
(p.~136) The notion of a non-standard model is complicated and highly
inconstructive [sic], but its properties, though surprising, are
rather simple and can therefore be understood by people who know very
little about model theory;
\item
(ibid.) Once you had shown by the paradigm of the calculus how it can
be used, many other applications were found by yourself as well as by
many other mathematicians.
\item
(p.~137) Model theory was in our time one of the abstract theories
which, so to say, floated in the air.  In a sense your work can be
considered as a return from this abstraction to {\bf concrete
applications} [emphasis added--authors].
\item
(ibid.) The general non-constructive theory of non-stand\-ard models
links its applications together into a harmoni[ous] whole.
\end{enumerate}
Heyting's final comment is particularly revealing.  He recognizes
that, while the model may be {\em non-constructive\/}, its
applications may indeed be {\em concrete\/} (in the terminology of
Subsection~\ref{postlem}, such applications possess post-LEM numerical
meaning).  Such a thesis is consonant with what we called {\em
numerical constructivism\/} in Subsection~\ref{51}.

In Heyting's case, an acute awareness of the shortcomings of
LEM-dependent mathematical investigations went hand-in-hand with an
appreciation of the accomplishments of classical mathematics in
general, and Robinson's mathematics in particular.  Such an awareness
is therefore not at odds with a recognition of the coherence of
post-LEM numerical meaning.  It is at odds, however, with an anti-LEM
radical constructivist stance, as enunciated by E.~Bishop.

\section{Constructivism, physics, and the real world}
\label{fourteen}

In this section we deal with challenges to anti-LEM constructivism
stemming from natural science.

\subsection{Bishop's view}
\label{81}

Bishop \cite[p.~4]{Bi85} speaks in favor of a goal of seeking
\begin{quote}
common ground in the researches of pure mathematics, applied
mathematics, and such mathematically oriented disciplines as {\bf
physics} [\dots] [emphasis added--authors].
\end{quote}
To what extent is such a commendable goal served well by
constructivism?

In the previous section, we mentioned the constructivist's problem of
accounting for applications of LEM-dominated mathematics.  This point
seems to have bothered the interviewer of Richman's {\em Interview\/}:
\begin{quote}
Here is a question that has bothered me ever since the first time I
read Bishop.  In what sense can your [constructive] real line be used
as a model for either space or time in physics? \cite[p.~252]{Ri}
\end{quote}
Richman replies as follows:
\begin{quote}
I don't know exactly how to respond to the [\dots] question, given that
I have already rejected the idea that constructivists have a real line
that is different from the classical real line.
\end{quote}

It needs to be understood what Richman meant by this cryptic remark.
For convenience, we will discuss the circle~$\R/\Z$ instead of~$\R$
itself.  The circle~$\R/\Z$ is constructively diffeomorphic (via the
exponential map) to the unit circle in the plane.  A classical
mathematician believes that the circle can be decomposed as the
disjoint union of the lower halfcircle (with one endpoint included,
and the other excluded) and its antipodal image (upper halfcircle).
Meanwhile, the constructivist believes that it is impossible to
decompose the circle as the disjoint union of a pair of antipodal
sets.  Aren't these different lines?  In a similar vein, a typical
physicist possesses a simple faith in the existence of infinitesimals
in our familiar number system, and Komkov \cite{Ko} explicitly
endorses Robinson's infinitesimals.  As we have documented in the
present text, Bishop and his followers are less eager to accept
Robinson's infinitesimals as part of our number system.  Aren't these
different ``real lines''?

At any rate, Richman's reply above%
\footnote{Now an anti-LEM constructivist may be satisfied with a
dodge, but the rest of us may wonder why Richman's 24-page {\em
Interview\/} has no bibliography, particularly in view of its
publication in a journal that describes itself as ``an organ for
rapid, low-cost, communication between historians of logic and
research logicians''.  Did history start with E.~Bishop?}
does not appear to deal with the issue of scientific applications as a
litmus test of the relevance of mathematical production.  Such an
issue is related to the {\em Quine-Putnam indispensability thesis\/},
see Feferman \cite[Section IIB]{Fe00}, M.~Colyvan~\cite{Co}, or
G.~Hellman~\cite{Hel}, see Subsection~\ref{fifteen} below.

In fact, Bishop himself was challenged on the relationship to physics
by G. Mackey, in the following terms:
\begin{quote}
Consider the foundational question in physics: what is the real
mathematics that the physicists are doing?  \cite[p.~515]{Bi75}
\end{quote}
In his reply, Bishop falls back on the issue of mathematical rigor:
\begin{quote}
physicists have told me that the sort of meaning that is appropriate
to physics is {\em not\/} to ask whether the mathematics in question
is rigorous.  Rather, it involves the relation of the results to the
{\bf real world} [emphasis added--authors].
\end{quote}

At first glance, Bishop's comment may seem to amount to an
unobjectionable observation concerning a higher standard of {\em
rigor\/} in mathematics.  However, what Bishop is emphasizing in his
essay is a notion of {\em meaning\/} which amounts to, as Hellman puts
it, {\em a philosophical principle of cognitive significance\/}
\cite[p.~239]{Hel93}.  Such a notion must, in Bishop's view
\cite[p.~509]{Bi75}, supercede and precede a notion of {\em truth\/},
as well as the related one of {\em rigor\/}.

The kind of (unabashedly LEM-dominated) mathematics that a physicist
practices is, as Bishop seems to admit matter-of-factly, successful in
relating to {\em the real world\/}, and therefore apparently {\em
meaningful\/}.  Such a situation apparently creates a tension with
Bishop's attempted identification of {\em meaning\/} with {\em
numerical meaning\/}, that can only be resolved in the context of a
numerical constructivism, a companion to classical mathematics,
recognizing the coherence of post-LEM numerical meaning (see
Subsection~\ref{postlem}).

In his {\em Constructivist manifesto\/}, Bishop wrote that
\begin{quote}
Weyl, a great mathematician who in practice suppressed his
constructivist convictions, expressed the opinion that idealistic
[i.e.~classical] mathematics finds its justification in its
applications in physics \cite[p.~10]{Bi67}.
\end{quote}
Bishop does not elaborate as to why he feels H.~Weyl was off the mark
here; this point will be discussed in Subsection~\ref{fifteen}.

\subsection{Hellman on indispensability}
\label{fifteen}

The challenge to constructivism stemming from the Quine-Putnam
indispesability thesis has been extensively pursued by G.~Hellman.
Hellman expresses an appreciation of the foundational significance of
constructivism in the following terms:
\begin{quote}
A turning point [\dots] came in the 1960s with the appearance of
Erret[t] Bishop's {\em Foundations of constructive analysis\/} (1967)
which succeeded in constructivizing large portions of classical
functional analysis of just the sort used in the physical sciences
(including the theory of operators on Banach and Hilbert spaces, and
even measure theory) \cite[p.~457]{Hel}.
\end{quote}
As to the question whether Bishop's constructivism can ultimately
serve the needs of the natural sciences, Hellman analyzes what he sees
as limitations of constructivism in scientific applications.  Hellman
\cite[p.~234]{Hel93a} invokes Plato's parable of the cave, apparently
reflecting his global perception of constructivism as an intellectual
house of cards.  However, one can recognize the coherence of the
intuitionist critique of the foundations of mathematics, while
rejecting the intellectual underpinnings of its insurrectional
narrative, in its mutually contradictory (realist and anti-realist)
versions, as discussed in Subsection~\ref{insurrection}.

In the context of general relativity theory, Hellman argues that, very
likely, the Hawking--Penrose singularity theorems for general
relativistic spacetimes are essentially non-constructive (Hellman in
\cite{He98}), see also (Billinge \cite{Bil}), see
Subsection~\ref{HPT}.

The philosopher of mathematics M.~Dummett rejects outright the idea
that ``empirical discoveries'' or applications should have any bearing
on the outcome of the debate opposing constructivism and classical
mathematics:
\begin{quote}
Intuitionists are engaged in the wholesale reconstruction of
mathematics, not to accord with empirical discoveries, not to obtain
more fruitful applications, but solely on the basis of philosophical
views concerning what mathematical statements are about and what they
mean [\dots] intuitionism will never succeed in its {\bf battle}
against rival, and more widely accepted, forms of mathematics unless
it can win the philosophical battle \cite[p.~viii]{Du} [emphasis
added--authors].
\end{quote}
The {\em battle\/} imagery is typical of the anti-realist type of
insurrectional narrative, already discussed in
Subsection~\ref{insurrection}.

Meanwhile, the philosopher of mathematics G.~Hellman holds, just as
categorically, that such a {\em philosophy-first\/} attempt is
``doomed''.%
\footnote{See footnote~\ref{doom} for Hellman's view.}
Awaiting the outcome of the philosophical battle, we note that
Dummett's stated intuitionist position qualifies for the anti-LEM
species.

More generally, variational problems tend to be resistant to efforts
at constructivizing, a point apparently acknowledged by Beeson when he
writes that
\begin{quote}
Calculus of variations is a vast and important field which lies right
on the frontier between constructive and non-constructive mathematics
\cite[p.~22]{Bee}.
\end{quote}
The problem is that the extreme value theorem is not available in the
absence of the law of excluded middle.  The extreme value theorem is
at the foundation of the calculus of variations.

As concrete examples, one could mention general existence results for
geodesics, minimal surfaces, constant curvature mean surfaces, and
therefore soap films and soap bubbles.  It seems implausible to argue
that existence results for soap films (referred to in the literature
as Plateau's problem) are meaningless merely because they lack
Bishop-style numerical meaning.

\begin{figure}
\includegraphics[height=2in]{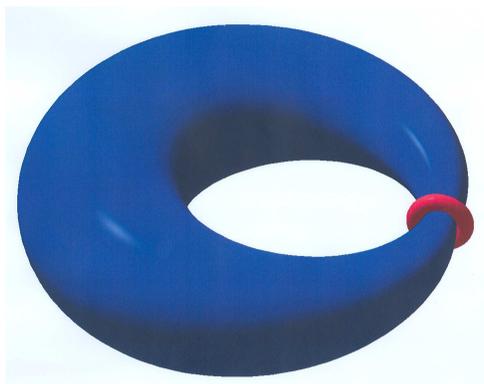}
\caption{A shortest loop on a torus}
\label{micro}
\end{figure}

\begin{question}
As a geometer, this writer would be interested in clarifying the
constructive status of the following phenomenon related to Loewner's
torus theorem of Example~\ref{loewner} in Subsection~\ref{postlem}.
Consider a rubber band ``clasping'' a torus in~$3$-space, and let it
shrink to a shortest position (see Figure~\ref{micro}).  Classically,
one has a description of the resulting object in terms of a minimizing
(possibly non-unique) closed geodesic.  Providing a constructive
description is not immediate, as it depends on general results of the
calculus of variations, as already described above.  To what extent
can the existence of a loop of least length of Figure~\ref{micro} be
recaptured constructively?  Pourciau's {\em Education\/}
\cite[p.~721]{Po99} advocates us to ``move from the simple to the less
simple'', but how plausible is it that such a loop could admit a {\em
simpler\/} description in terms of a statement about the {\em
integers\/}, as advocated by {\em Integrity\/} (see glossary in
Subection~\ref{glossary})?
\end{question}

\subsection{Novikov's perspective}
\label{83}

S. P. Novikov's essay \cite{No2}, in the ``crisis" tradition of Weyl
\cite{We21} and Bishop \cite{Bi75}, deals with the crisis in the
cooperation of physicists and mathematicians, in Russia and in the
West.  What Novikov makes excruciatingly clear is the disastrous
effect of such a divorce on mathematics itself.  The crucial role of
variational principles is made clear already on the first page of
Novikov's essay.  A broad range of subjects illustrating the symbiotic
relationship between the two fields can be found in Grattan-Guinness
\cite{Gra08}.

Grattan-Guinness emphasizes the role of {\em analogies drawn from
other theories\/}, in theory-building in both mathematics and the
sciences.  From his viewpoint, the physicists that Novikov has in mind
have been going by the book (Grattan-Guinness's book, at any rate),
applying mathematical machinery {\em by analogy\/}, outside its
original domain of applicability, and doing so with great success (in
terms of predicting experiments, etc.); yet no theoretical
justification for their calculations, in terms of traditional
mathematics, has followed.  What is worse, the gap is growing wider,
according to Novikov.  Traditional mathematical foundations, whether
classical or intuitionist, are proving to be inadequate for the job of
accounting for the progress in physics.  Accordingly, Novikov is
critical of set-theoretic foundations, going as far as criticizing
Kolmogorov himself, for systematic efforts to introduce a
set-theoretic approach in secondary education.  Foundational
insurrections \`a la Brouwer, Bishop, and Dummett (see
Subsection~\ref{insurrection}) may amount to irrelevant distractions,
from Novikov's viewpoint.

\subsection{The Hawking--Penrose theorem}
\label{HPT}

The Hawking--Penrose theorem in relativity theory has been the subject
of something of a controversy in its own right, see G.~Hellman
\cite{He98}, H.~Billinge \cite{Bil}, and E.~Davies \cite{Da05}.
Hellman argued that the Hawking--Penrose theorem is an important result
which very likely does not have a constructive analog.  As we explain
below, Billinge's criticism of Hellman's argument is based on
conflating two separate issues (geodesic incompleteness, on the one
hand, and a hypothetical description of a ``singularity'', on the
other), due to a mathematical error on her part.

Davies \cite[p.~268]{Da05} mentions the following points:
\begin{enumerate}
\item
\label{d1}
``Hellman \cite{He98} [\dots] showed that the Hawking--Penrose
singularity theorem is not constructively valid.''
\item
\label{d2}
``the theorem has been extremely influential in the subject.''
\item
\label{d3}
``It is embarrasssing for a radical constructivist [\dots]''
\item
\label{d4}
``It remains extremely hard to say much about the nature of the
singularities [\dots]  It is very likely that if a detailed description
of the singularities becomes possible classically, that description
will also be constructively valid.''
\end{enumerate}

We will return to Davies' remarks after reviewing the relevant
mathematical details.

The Hawking--Penrose singularity theorem can be thought of as a
semi-Riemannian analog of the Myers theorem \cite{My}.  The latter is
a result in Riemannian geometry.  The result in question is a bound,
modulo a suitable geometric hypothesis on the Ricci curvature, on the
distance from a point to its conjugate locus.  The traditional
formulation of the Myers theorem relies upon the concept of a geodesic
and a Jacobi field, both of which are solutions of variational
problems.  Myers found a bound on the distance along a
geodesic~$\gamma$, to the nearest vanishing point, called a {\em
conjugate point\/}, of a Jacobi field along~$\gamma$.

Returning to the Lorenzian case, note that Penrose does not even use
the term {\em singularity\/} in formulating the result that came to be
known as the Hawking--Penrose theorem.  Rather, what one claims is the
existence of a ``past-endless geodesic [\dots] which has finite
length'' \cite[p.~69]{Pe}.  The proof is indirect, i.e.~a proof by
contradiction.  Penrose describes one of the main ingredients in the
proof as the Raychaudhuri effect.  The latter is discussed in
\cite[item 7.21, p.~63]{Pe}.  Penrose points out that 
\begin{quote}
[f]or manifolds with a positive definite metric, essentially the same
effect had been studied earlier by Myers (Penrose
\cite[p.~64]{Pe}).
\end{quote}

Thus, the formulation of the Hawking--Penrose theorem in terms of
``singularities'' is superfluous to its purely mathematical content,
namely the fact that space-time must be (time-like or null)
geodesically incomplete.

Note that there exist examples of {\em compact\/} Lorentz manifolds
which are geodesically incomplete \cite{HE}.  Penrose places {\em
singularities\/} in quotation marks when he points out that
\begin{quote}
the physical implication of the theorem is that ``singularities''
(i.e., causal geodesic incompleteness) would be expected to arise
whenever such a collapse takes place \cite[p.~71]{Pe}.
\end{quote} 
(Note the significant {\em i.e.\/}) Whether or not geodesic
incompleteness of space-time can be understood in terms of a suitable
singular limiting object, boundary, or hole, is a separate question
(extensively discussed in the physics literature, see references cited
in \cite[p.~213]{Wal}), not germane to the problem of the foundational
status (classical or intuitionist) of the Hawking--Penrose theorem.

Furthermore, while a suitable bound on the Ricci curvature does allow
one to prove the existence of limiting objects in a Riemannian context
(see Gromov \cite[Theorem~5.3, p.~275]{Gr4}), such limiting objects
are known to be arbitrarily pathological in general, unless one
imposes additional geometric hypotheses (as in the Anderson-Cheeger
theory~\cite{AC}, \cite{Ta}), or, alternatively, leaves the domain of
differential geometry and imposes conditions of algebraicity (which do
not appear to be justified physically) on the manifolds being studied.

Foundationally or otherwise, finding fault with the Hawking--Penrose
theorem for failing to describe a ``singularity'', is akin to finding
fault with the Big Bang theory for failing to describe ``what came
before the bang''.

Thus, Billinge's objection to Hellman's text \cite{He98}, on the
grounds of a lack of a constructive description of the singularity
(echoed in item~\eqref{d4} above), does not appear to be based on a
detailed understanding of the mathematics.  The following phrase
appears in Billinge:
\begin{quote}
a spacetime is singular if it is timelike or null geodesically
incomplete, {\bf that is,} if it has at least one timelike or null
geodesic that has a point at which it comes to an end
\cite[p.~308]{Bil} [emphasis added--authors].
\end{quote}

The words ``that is'' that we have boldfaced, apparently presented by
Billinge as a mathematical definition, constitute an elementary
mathematical error.  The spacetime in question does {\em not\/}
contain geodesics that ``come to an end'', at any point of the
spacetime.  Rather, the geodesics are, as Penrose puts it in
\cite[p.~69]{Pe}, {\em past-endless\/} of finite length, i.e.~they run
off to (negative) infinity in finite time (in the Lorentzian case,
there is no analogue of the Hopf-Rinow theorem).  Naturally, if one
thinks of geodesics as running into some kind of a singular
``endpoint'', or barrier, in spacetime, one can easily develop an
impression that the relation of geodesic incompleteness to
``singularities'' is far more immediate than it really is.

Billinge's main source for relativity is the textbook by
Wald~\cite{Wal}.  The textbook contains the following remarks, which
express the distinction between {\em singularity\/} and {\em geodesic
incompleteness\/} clearly:
\begin{quote}
it is extremely difficult to give a satisfactory general notion of [a
`singularity'].  We provide motivation for the notion of timelike and
null geodesic incompleteness as a criterion for the presence of a
singularity [\dots] It is this criterion which is used in the
singularity theorems \cite[p.~212]{Wal}.
\end{quote}
Furthermore, on page 214, Wald clearly states that
\begin{quote}
[u]ntil a satisfactory definition can be produced, we must abandon the
notion of a singularity as a ``place''.
\end{quote}
Describing a singularity as a ``place'' is, of course, precisely what
Billinge did at \cite[p.~308]{Bil}.  Such an error could have been
written off as a mere slip, were it not for the further evidence of
insufficient background presented above.

It should be mentioned that, at variance with Davies' item~\eqref{d1}
above, Hellman never purported to {\em prove\/} that ``the
Hawking--Penrose singularity theorem is not constructively valid'',
but only that he finds this ``very likely'', based on the
non-constructive nature of the proofs found in the literature, and the
general difficulties of accounting for the calculus of variations
constructively (in this case, accounting for geodesics and Jacobi
fields), as already discussed in Subsection~\ref{fifteen}.

N.~Tennant, writing on {\em Logic, mathematics, and the natural
sciences\/}, claims that
\begin{quote}
a constructivist version of a mathematical theory is adequate for all
the applications to be made of the theory within natural science
\cite[p.~1145]{Te06}.
\end{quote}
Tennant claims that such a constructivist version is adequate for
applications.  Or perhaps it isn't?  Tennant's claim apparently flies
in the face of everything we have written here.  In fact, Tennant's
claim is not substantiated in his article.  Tennant continues by
explaining that
\begin{quote}
to be able to produce all possible {\bf refutations} of empirical
theories, the underlying logic can [\dots] be taken very weak
\cite[p.~1146]{Te06} [emphasis added--authors].
\end{quote}
Furthermore he concedes:
\begin{quote}
Nor do we intend to say anything about confirmation or
probabilification of hypotheses by evidence \cite[p.~1149]{Te06}.
\end{quote}
The evidence furnished by astrophysical observation tends to confirm
the hypotheses of the Hawking--Penrose theorem, but Tennant's framework
does not deal with this type of scientific insight.

It is possible that, should further empirical observations one day
undermine the premises of the Hawking--Penrose theorem, the resulting
{\em refutation\/} may turn out to be adequately expressible in
Tennant's version of intuitionistic logic (this does not quite follow
from \cite{Te06} as Tennant limits the discussion to first order
logic).  However, this does not alter the fact that expressing the
Hawking--Penrose theorem (let alone proving it) remains an elusive
goal, intuitionistically speaking, with the attendant loss of
scientific insight provided by the theorem.

The Hawking--Penrose theorem, described by Davies as {\em extremely
influential\/} (see item~\eqref{d2} above), remains a challenge to
Bishopian constructivism (see item~\eqref{d3} above), as do the
numerous variational principles mentioned by Novikov \cite{No2}, see
also Grattan-Guinness \cite{Gra08}.

\section{Epilogue}

In 1973, Robinson delivered the {\em Brouwer memorial lecture\/}
\cite{Ro73} on the subject of {\em Standard and nonstandard number
systems\/}.%
\footnote{See Heyting's introduction, quoted at length in
Subsection~\ref{APR}.}
J. Dauben quotes Robinson as follows:
\begin{quote}
Brouwer's intuitionism is closely related to his conception of
mathematics as a {\bf dynamic activity} of the human intellect [\dots]
This is a conception for which I have some sympathy and which, I
believe, is acceptable to many mathematicians who are not
intuitionists \cite[p.~461]{Da95} [emphasis added--authors].
\end{quote}
Robinson concludes that
\begin{quote}
the dynamic evolution of mathematics is an ongoing process not only at
the summit [\dots] but also at the more basic level of our number
systems, which [tomorrow] may seem as eternal to a new generation as
yesterday's technological innovations are in the eyes of a child of
today.
\end{quote}

Have today's constructivists been faithful to a conception of
mathematics as a {\em dynamic activity\/}, as Robinson puts it?  We
would like to suggest that the possibility of bridging the gap between
constructivism and non-standard analysis (as well as, indeed, the rest
of classical mathematics) can be analyzed in the context of the
dichotomy of a numerical constructivism {\em versus\/} an anti-LEM
constructivism (see Section~\ref{51}).  Namely, only the former,
consistent with Heyting's brand of intuitionism, seems sufficiently
{\em dynamic\/} to answer the hopes~\cite{Sc} of {\em reuniting the
antipodes\/}.

\appendix

\section{Rival continua}
\label{rival}

\begin{figure}
\[
\xymatrix@C=95pt{{} \ar@{-}[rr] \ar@{-}@<-0.5pt>[rr]
\ar@{-}@<0.5pt>[rr] & {} \ar@{->}[d]^{\hbox{st}} & \hbox{\quad
B-continuum} \\ {} \ar@{-}[rr] & {} & \hbox{\quad A-continuum} }
\]
\caption{\textsf{Thick-to-thin: taking standard part (the thickness of
the top line is merely conventional)}}
\label{31}
\end{figure}

A Leibnizian definition of the derivative as the infinitesimal
quotient
\[
\frac{\Delta y}{\Delta x},
\] 
whose logical weakness was criticized by Berkeley, was modified by
A.~Robinson by exploiting a map called {\em the standard part\/},
denoted~``st'', from the finite part of a B-continuum (for
``Bernoullian''), to the A-continuum (for ``Archimedean''), as
illustrated in Figure~\ref{31}.

We will denote such a B-continuum by the new symbol \RRRhuge. We will
also denote its finite part, by
\[
\RRR_{<\infty} = \left\{ x\in \RRR : \; |x|<\infty \right\}, 
\]
so that we have a disjoint union~$\RRR=\RRR_{<\infty}\cup\RRR_\infty$,
where~$\RRR_\infty$ denotes the set of inverses of nonzero
infinitesimals.

The map ``st'' sends each finite point~$x\in \RRR$, to the real point
st$(x)\in \R$ infinitely close to~$x$:
\[
\xymatrix{\quad \RRR_{{<\infty}}^{~} \ar[d]^{{\rm st}} \\ \R}
\]
Robinson's answer to Berkeley's {\em logical criticism\/}%
\footnote{See footnote~\ref{sherry}.}
is to define the derivative as
\[
\hbox{st} \left( \frac{\Delta y}{\Delta x} \right),
\]
instead of~$\Delta y/\Delta x$ as in Leibniz.  ``However, this is a
small price to pay for the removal of an inconsistency''
\cite[p~266]{Ro66}.

\begin{figure}
\includegraphics[height=2in]{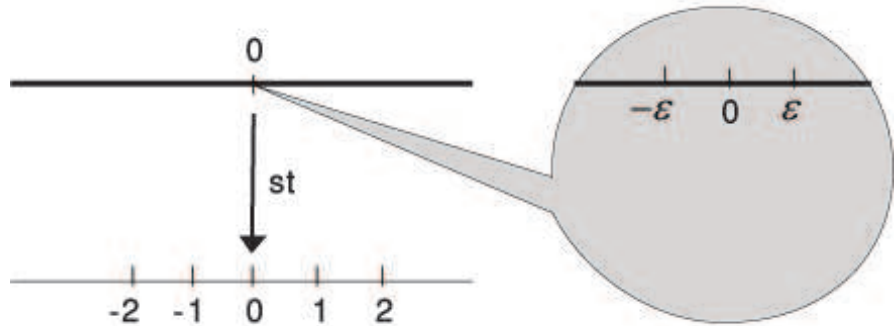}
\caption{Zooming in on infinitesimal~$\epsilon$}
\label{tamar}
\end{figure}

We illustrate the construction by means of an infinite-resolution
microscope in Figure~\ref{tamar}.

Note that both the term ``hyper-real field'', and an
ultrapower construction thereof, are due to E.~Hewitt in 1948, see
\cite[p.~74]{Hew}.  The transfer principle allowing one to extend
every first-order real statement to the hyperreals, is due to
J.~{\L}o{\'s} in 1955, see \cite{Lo}.  Thus, the Hewitt-{\L}o{\'s}
framework allows one to work in a B-continuum satisfying the transfer
principle.  To elaborate on the ultrapower construction of the
hyperreals, let~$\Q^\N$ denote the space of sequences of rational
numbers.  Let~$\left( \Q^\N \right)_C$ denote the subspace consisting
of Cauchy sequences.  The reals are by definition the quotient
field
\[
\R:= \left. \left( \Q^\N \right)_C \right/ \mathcal{F}_{\!n\!u\!l\!l},
\]
where the ideal~$\mathcal{F}_{\!n\!u\!l\!l}$ contains all the null
sequences.  Meanwhile, an infinitesimal-enriched field extension
of~$\Q$ may be obtained by forming the quotient
\[
\left.  \Q^\N \right/ \mathcal{F}_{u},
\]
See Figure~\ref{helpful}.  Here a sequence~$\langle u_n \rangle$ is
in~$\mathcal{F}_{u}$ if and only if the set
\[
\{ n \in \N : u_n = 0 \}
\]
is a member of a fixed ultrafilter.%
\footnote{An ultrafilter on~$\N$ can be thought of as a way of making
a systematic choice, between each pair of complementary infinite
subsets of~$\N$, so as to prescribe which one is ``dominant'' and
which one is ``negligible''.  Such choices have to be made in a
coherent manner, e.g., if a subset~$A\subset \N$ is negligible then
any subset of~$A$ is negligible, as well.  The existence of
ultrafilters was proved by Tarski \cite{Tar}, see Keisler
\cite[Theorem~2.2]{Ke08}.}
To give an example, the sequence~$\left\langle \tfrac{(-1)^n}{n}
\right\rangle$ represents a nonzero infinitesimal, whose sign depends
on whether or not the set~$2\N$ is a member of the ultrafilter.  To
obtain a full hyperreal field, we replace~$\Q$ by~$\R$ in the
construction, and form a similar quotient
\[
\RRR:= \left.  \R^\N \right/ \mathcal{F}_{u}.
\]

A more detailed discussion of the ultrapower construction can be found
in M.~Davis~\cite{Da77}, Gordon {\em et al.\/} \cite{GKK}, and
Giordano and Katz \cite{GK11}.  See also B\l aszczyk~\cite{Bl} for
some philosophical implications.  More advanced properties of the
hyperreals such as saturation were proved later, see Keisler
\cite{Kei} for a historical outline.  A helpful ``semicolon'' notation
for presenting an extended decimal expansion of a hyperreal was
described by A.~H.~Lightstone~\cite{Li}.  A discussion of
infinitesimal optics is in K.~Stroyan \cite{Str},
H.~J.~Keisler~\cite{Ke}, D.~Tall~\cite{Ta80}, and L.~Magnani and
R.~Dossena~\cite{MD, DM}.  P.~Ehrlich recently constructed an
isomorphism of maximal surreals and hyperreals \cite{Eh11}.
Applications of the B-continuum range from aid in teaching calculus
\cite{El, KK10a, KK10b, Ta91, Ta09a} to the Bolzmann equation (see
L.~Arkeryd~\cite{Ar81, Ar05}); modeling of timed systems in computer
science (see H.~Rust \cite{Rust}); mathematical economics (see
Anderson~\cite{An00}); mathematical physics (see Albeverio {\em et
al.\/} \cite{Alb}); etc.

\begin{figure}
\[
\xymatrix{ && \left( \left. \Q^\N \right/ \mathcal{F}_{\!u}
\right)_{<\infty} \ar@{^{(}->} [rr]^{} \ar@{->>}[d]^{\rm st} &&
\RRR_{<\infty} \ar@{->>}[d]^{\rm st} \\ \Q \ar[rr] \ar@{^{(}->} [urr]
&& \R \ar[rr]^{\simeq} && \R }
\]
\caption{\textsf{An intermediate field~$\left. \Q^\N \right/
\mathcal{F}_{\!u}$ is built directly out of~$\Q$}}
\label{helpful}
\end{figure}

In 1990, Hewitt reminisced about his 
\begin{quote}
efforts to understand the ring of all real-valued continuous [not
necessarily bounded] functions on a completely regular~$T_0$-space. I
was guided in part by a casual remark made by Gel'fand and Kolmogorov
(Doklady Akad. Nauk SSSR 22 [1939], 11-15). Along the way I found a
novel class of real-closed fields that superficially resemble the real
number field and have since become the building blocks for nonstandard
analysis. I had no luck in talking to Artin about these hyperreal
fields, though he had done interesting work on real-closed fields in
the 1920s. (My published ``proof" that hyperreal fields are
real-closed is false: John Isbell earned my gratitude by giving a
correct proof some years later.)  [\dots] My ultra-filters also struck
no responsive chords. Only Irving Kaplansky seemed to think my ideas
had merit.  My first paper on the subject was published only in 1948
\cite{He90}
\end{quote}
(Hewitt goes on to detail the eventual success and influence of his
1948 text).  Here Hewitt is referring to Isbell's 1954 paper
\cite{Isb}, proving that Hewitt's hyper-real fields are real closed.
Note that a year later, \Los~\cite{Lo} proved the general transfer
principle for such fields, implying in particular the property of
being real closed, the latter being first-order.

\section*{Acknowledgments}

We are grateful to Geoffrey Hellman for insightful comments that
helped improve an earlier version of the manuscript, to Douglas
Bridges for a number of helpful remarks, and to Claude LeBrun for
expert comments on the Hawking--Penrose theorem.  Hilton Kramer's
influence is obvious throughout.


\vfill\eject

\medskip\noindent
{\bf Karin Usadi Katz} has taught mathematics at Michlelet Banot
Lustig, Ramat Gan, Israel.  Two of her joint studies with Mikhail Katz
were published in {\em Foundations of Science\/}: ``A Burgessian
critique of nominalistic tendencies in contemporary mathematics and
its historiography" and ``Stevin numbers and reality", online
respectively at

http://dx.doi.org/10.1007/s10699-011-9223-1 and at

http://dx.doi.org/10.1007/s10699-011-9228-9

A joint study with Mikhail Katz entitled ``Cauchy's continuum'' is due
to appear in {\em Perspectives on Science\/} 2011, see
http://arxiv.org/abs/1108.4201

\medskip\noindent
{\bf Mikhail G. Katz} is Professor of Mathematics at Bar Ilan
University, Ramat Gan, Israel.  Two of his joint studies with Karin
Katz were published in {\em Foundations of Science\/}: ``A Burgessian
critique of nominalistic tendencies in contemporary mathematics and
its historiography" and ``Stevin numbers and reality", online
respectively at

http://dx.doi.org/10.1007/s10699-011-9223-1 and at

http://dx.doi.org/10.1007/s10699-011-9228-9

A joint study with Karin Katz entitled ``Cauchy's continuum'' is due
to appear in {\em Perspectives on Science\/}, 2011, see
http://arxiv.org/abs/1108.4201

A joint study with Alexandre Borovik entitled ``Who gave you the
Cauchy--Weierstrass tale?  The dual history of rigorous calculus''
appeared in {\em Foundations of Science\/}, online at

http://dx.doi.org/10.1007/s10699-011-9235-x

A joint study with David Tall, entitled ``The tension between
intuitive infinitesimals and formal mathematical analysis'', is due to
appear as a chapter in a book edited by Bharath Sriraman, see

\noindent
http://www.infoagepub.com/products/Crossroads-in-the-History-of-Mathematics

\end{document}